\documentclass[a4paper]{article}
\usepackage{amsmath,amssymb,fancybox,longtable,amscd,multirow,theorem}
\usepackage{caption,setspace,framed,arydshln,graphicx}
\newtheorem{thm}{\textsc{Theorem}}[section]

\newtheorem{lem}{\textsc{Lemma}}[subsection]

\def\proof{\textsc{Proof. }}
\def\QED{$\Box$} 

\def\wt {\mathop{\mathrm{wt}}\nolimits}

\theorembodyfont{\rmfamily}
\newtheorem{remark}{\textsc{Remark}}

\theorembodyfont{\rmfamily}

\begin{document} 
\title{Polytope duality for families of $K3$ surfaces and coupling}
\author{Makiko Mase}
\date{\empty}
\maketitle
\abstract{We study a relation between coupling introduced in \cite{Ebeling} and the polytope duality  among families of $K3$ surfaces.}


\section{Introduction}
A notion of coupling is introduced by Ebeling~\cite{Ebeling} as a tone-down of the duality of weight systems by Kobayashi~\cite{Kobayashi}. 
It is proved that the duality is also ``polar dual'', in the sense that certain rational polytopes associated to weight systems are dual.  
In~\cite{Ebeling}, there is given a list of coupling pairs for 95 weight systems of simple $K3$ hypersurface singularities classified by Yonemura~\cite{Yonemura}, and it is proved that the duality induces Saito's duality, which is a relation between the zeta functions of the Milnor fibre of the singularities. 
It is interesting to note that these dualities translate a famous mirror symmetry. 

Instead of rational polytopes, we are interested in Batyrev's toric mirror symmetry~\cite{BatyrevMirror} for integral polytopes in this article. 
As a generalisation of Arnold's strange duality for unimodal singularities, Ebeling and Takahashi~\cite{EbelingTakahashi11} defined a notion of strange duality for invertible polynomials. 
It is studied by Mase and Ueda~\cite{MU} that the strange duality for bimodal singularities defined by invertible polynomials extends to the polytope duality among families of $K3$ surfaces. 

The polytope duality is focusing on more details in geometry of $K3$ surfaces such as resolution of singularities, as a compactification of some singularities in three dimensional space which should affect the geometry of the surfaces while the polar duality in~\cite{Ebeling,Kobayashi} is determined only by the weight systems and thus, it is quite global. 
There are some profiles in coupling as is explained in~\cite{Ebeling} from the viewpoint of mirror symmetry, in particular, in terms of the Milnor fibres. 
It is expected that a study of the polar duality associated to coupling gives another explanation to coupling by Batyrev's mirror symmetry. 
In turn, we expect some relation can be extracted between the Milnor fibres of singularities and the geometry of associated $K3$ surfaces. 

Motivated by this, and focusing on coupling, we consider the following problem. \\

\noindent
{\sc Problem. } 
Let $a$ and $b$ be weight systems that are coupling pair and their weighted magic square is given by polynomials $f$ and $f'$, respectively. 
Determine whether or not there exist reflexive polytopes $\Delta$ and $\Delta'$, and projectivisations $F$ and $F'$ of $f$ and $f'$ in the weighted projective spaces $\mathbb{P}(a)$ and $\mathbb{P}(b)$, respectively, such that they are {\it polytope dual} in the sense that they satisfy the following conditions:
\[
\Delta^*\simeq \Delta',\quad \Delta_F\subset\Delta\subset\Delta_{a},\quad \Delta_{F'}\subset\Delta'\subset\Delta_{b}. 
\]
Here $\Delta_{F}$ and $\Delta_{F'}$ are Newton polytopes of $F$ and $F'$, respectively, and $\Delta_{a}$ and $\Delta_{b}$ are polytopes that define the weighted projective spaces $\mathbb{P}(a)$ and $\mathbb{P}(b)$. 

The main theorem of the article, which is proved in Section 5 is stated:

\noindent
{\bf Theorem \ref{MainThm}}\quad 
{\it 
Any coupling pairs in Yonemura's list extend to the polytope dual except the following three  pairs of weight systems : 
$(1,3,4,7;15)$ (self-coupling),\, $(1,3,4,4;12)$ (self-coupling),\, and $(1,1,3,5;10)$ and $ (3,5,11,19;38)$. 
An explicit choice of reflexive polytopes is given in Table~\ref{ListMainThm}. 
} \\

In section 2, we recall the definitions concerning the weighted projective spaces and the strange duality. 
In section 3, we recall the definition of coupling. 
In section 4, we explain the polytope duality after recalling necessary notions of toric geometry. 

\begin{ackn}
\textnormal{
The author thanks to Professor Wolfgang Ebeling for his suggestion of the study and discussions. }
\end{ackn}

\section{Preliminary}
A {\it $K3$ surface} is a compact complex $2$-dimensional non-singular algebraic variety with trivial canonical bundle and irregularity zero. 

Let $(a_0,\,\ldots ,\, a_n)$ be a well-posed $(n+1)$-tuple of positive integers, that is, $a_0\leq\ldots\leq a_n$, and any $n$-tuples out of them are coprime. 
Recall that the {\it weighted projective space} $\mathbb{P}(a)=\mathbb{P}(a_0,\,\ldots ,\, a_n)$ with {\it weight} $a=(a_0,\,\ldots ,\, a_n)$ is defined by 
\[
\mathbb{P}(a)=\mathbb{P}(a_0,\,\ldots ,\, a_n) := \mathbb{C}^{n+2}\backslash\{ 0\} \biggr\slash \sim, 
\]
where $(x_0,\,\ldots ,\, x_n) \sim (y_0,\,\ldots ,\, y_n)$ if there exists a non-zero complex number $\lambda$ such that
\[
(y_0,\,\ldots ,\, y_n) = (\lambda^{a_0}x_0,\,\ldots ,\, \lambda^{a_n}x_n)
\]
holds. 
We call $a_i$ the {\it weight} of the variable $x_i$, and denote it by $\wt x_i$. 

In case $n=3$, we fix a system of variables $W,\, X,\, Y,\, Z$ of the weighted projective space $\mathbb{P}(a_0,\, a_1,\, a_2,\, a_3)$ with weights 
\[
\wt{W}=a_0,\quad
\wt{X}=a_1,\quad
\wt{Y}=a_2,\quad
\wt{Z}=a_3. 
\]
We say a polynomial $F$ in $\mathbb{P}(a_0,\, a_1,\, a_2,\, a_3)$ is an {\it anticanonical section} if $F$ is of degree $d:=a_0+a_1+a_2+a_3$. 
The tuple $(a_0,\, a_1,\, a_2,\, a_3;\, d)$ is called a {\it weight system}. 

By~\cite{DolgachevWP}, the anticanonical sheaf of $\mathbb{P}(a)$ is isomorphic to $\mathcal{O}_{\mathbb{P}(a)}({-}d)$. 
All weight systems that give simple $K3$ hypersurface singularities are classified by Yonemura~\cite{Yonemura}. 
Namely, if a weight $a$ is in Yonemura's list, general anticanonical sections of $\mathbb{P}(a)$ are birational to $K3$ surfaces. 
Thus, one can consider families of $K3$ surfaces. 


For a polynomial $f$ in three variables, a polynomial $F$ in the weighted projective space $\mathbb{P}(a)$ is called a {\it projectivisation of $f$} if there exists a linear form $l$ in $\mathbb{P}(a)$ such that 
\[
f = F|_{l=0}
\]
holds. 
In this case, the form $l$ is called a {\it section} of $f$ for $F$.  

\section{Coupling}
Recall the definition of {\it coupling} for weight systems with three entries. 
Let $w=(w_1,\,w_2,\,w_3;\, d)$ and $w'=(w'_1,\,w'_2,\,w'_3;\, h)$ be weight systems, with weights being well-posed. 

A {\it weighted magic square} $C$ for the weight systems $w$ and $w'$ is a square matrix of size $3$ that satisfies relations 
\[
C\,{}^t\!(w_1 \; w_2 \; w_3) = {}^t\!(d \; d \; d)
\textnormal{ and }
(w'_1 \; w'_2 \; w'_3)\,C = (h \; h\; h). 
\]
The pair of weight systems $(w,\, w')$ is called {\it coupled} if $C$ is {\it almost primitive}, that is, if $|{\det{C}}| = \left(d-\sum_{i=1}^3w_i\right)h=\left(h-\sum_{i=1}^3w'_i\right)d$ hold. 
The pair of weight systems $(w,\, w')$ is {\it strongly coupled} if it is coupled and the weighted magic square $C$ has entries zero in every column and row. 

Thus, one can assign polynomials $f$ and $f'$ to $C=(c_{ij})$ in such a way that
\[
f=\sum_{i=1}^3x^{c_{i1}}y^{c_{i2}}z^{c_{i3}}, 
\quad 
f'=\sum_{i=1}^3x^{c_{1i}}y^{c_{2i}}z^{c_{3i}}. 
\]
In other words, there exisis weight systems $(a_0,\, a_1,\, a_2,\, a_3;\, d)$ and $(b_0,\, b_1,\, b_2,\, b_3;\, h)$ such that there exist $i,\,j\in\{0,1,2,3\}$ with properties 
\[
|{\det{C}}| = ha_i = kb_j, 
\textnormal{ and }
{}^t\!A_{f'} = A_f. 
\]
Define an anticanonical section $F$ of weight system $(a_0,\, a_1,\, a_2,\, a_3; \, d)$ so that $l$ is the section of $f$ for $F$, 
where $l$ is a linear form defined by 
\[
l=
\begin{cases} 
w^{h/a_0} & \textnormal{if } |{\det{C}}|=ha_0, \\
x^{h/a_1} & \textnormal{if } |{\det{C}}|=ha_1,  \\
y^{h/a_2} & \textnormal{if } |{\det{C}}|=ha_2,  \\
z^{h/a_3} & \textnormal{if } |{\det{C}}|=ha_3. 
\end{cases}
\]
Note that the choice of variables is different from the original Ebeling's paper \cite{Ebeling}. 
And then define a polynomial $F'$ so that
\[
{}^t\!A_{F'} = A_F
\]
holds. 
Note that $F'$ is a projectivisation of $f'$, and an anticanonical section in the weighted projective space of weight system $(b_0,\, b_1,\, b_2,\, b_3;\, h)$. 

Ebeling~\cite{Ebeling}(Tables 2 and 3) gives (strongly) coupling pairs among weighted systems in Yonemura's list.

\section{Duality of polytopes}
Let $M$ be a lattice of rank $3$, and $N$ be its dual lattice $\mathrm{Hom}_{\mathbb{Z}}(M,\, \mathbb{Z})$ that is again of rank $3$. 
A {\it polytope} is a convex hull of finite number of points in $M\otimes\mathbb{R}$. 
If vertices of a polytope $\Delta$ are $v_1,\ldots,v_r$, we denot it by 
\[
\Delta = \mathrm{Conv}\{ v_1,\ldots,v_r\}. 
\]
We call a polytope {\it integral} if all the vertices of the polytope are in $M$. 
For a polytope $\Delta$, define the {\it polar dual polytope} $\Delta^*$ by 
\[
\Delta^* := \{y\in N\otimes\mathbb{R}\, |\, \langle y,x\rangle_{\mathbb{R}}\geq-1\, \textnormal{for all}\,  x\in\Delta\},
\]
where $\langle \, ,\, \rangle$ is a natural pairing $N\times M\to\mathbb{Z}$, and $\langle \, ,\, \rangle_{\mathbb{R}}$ is the extension to $\mathbb{R}$-coefficients. 
Let $\Delta$ be an integral polytope that contains the origin in its interior as the only lattice point. 
The polytope $\Delta$ is {\it reflexive} if the polar dual $\Delta^*$ is also an integral polytope. 

Recall an interesting property of reflexive polytopes related to $K3$ surfaces due to Batyrev~\cite{BatyrevMirror}: 
\begin{thm}\cite{BatyrevMirror}
Denote by $\mathbb{P}_{\!\Delta}$ the toric $3$-fold associated to an integral polytope $\Delta$. 
The following conditions are equivalent. 
\begin{enumerate}
\item[$(1)$] The polytope $\Delta$ is reflexive. 
\item[$(2)$] General anticanonical sections of $\mathbb{P}_{\!\Delta}$ are birational to $K3$ surfaces. \QED
\end{enumerate}
\end{thm}

In particular, the weighted projective space $\mathbb{P}(a)$ with weight system $a=(a_0,\, a_1,\, a_2,\, a_3)$ in Yonemura's list is a toric Fano $3$-fold determined by a reflexive polytope $\Delta^{(n)}$ in the $\mathbb{R}$-extension of the lattice 
\[
M_{n}:=\left\{ (i,\, j,\, k,\, l)\in\mathbb{Z}^4 \, |\, a_0i + a_1j + a_2k + a_3l=0\right\}, 
\]
where the weight system $a$ is assigned No. $n$ in Yonemura's list. 
The anticanonical sections are weighted homogeneous polynomial of degree $d:=a_0 + a_1 + a_2 + a_3$, thus, there is a one-to-one correspondence between a lattice point $(i,\, j,\, k,\, l)$ in $M_n$ and a rational monomial $W^{i+1}X^{j+1}Y^{k+1}Z^{l+1}$. 
In this way, once a $\mathbb{Z}$-basis is taken for $M_n$, we identify lattice points in $\Delta^{(n)}$ and monomials of weighted degree $d$. 

\section{Main Result}
In this section, we prove the main theorem. 

\begin{thm}\label{MainThm}
Any coupling pairs in Yonemura's list extend to the polytope dual except the following three  pairs of weight systems : 
$(1,3,4,7;15)$ (self-coupling),\, $(1,3,4,4;12)$ (self-coupling),\, and $(1,1,3,5;10)$ and $ (3,5,11,19;38)$. 
An explicit choice of reflexive polytopes is given in Table~\ref{ListMainThm}. 
\begingroup
\setlength{\arrayrulewidth}{.1pt}
\tiny
\begin{longtable}[htp]{p{1.8mm}|
p{30mm}@{\hspace{7mm}}
p{10mm}@{\hspace{7mm}}
p{10mm}@{\hspace{7mm}}
p{35mm}}
{\rm No.} & $\Delta'$ & $b;h$ & $a;d$ & $\Delta$ \\
\hline
$1.$ 
& \shortstack{\\$Z^2,\, W^{42},\, X^7,\, Y^3$} 
& $1,6,14,21;42$ 
& $1,6,14,21;42$  
&  \shortstack{\\$Z^2,\, W^{42},\,X^7,\, Y^3$}  \\
\hline
$2.$ 
&  \shortstack{\\$WZ^2,\, W^{21},\, X^7,\, Y^3$} 
& $1,3,7,10;21$ 
& $1,6,14,21;42$ 
&  \shortstack{\\$Z^2,\, W^{42},\,X^7,\, Y^3$}  \\
\hline
$3.$ 
&  \shortstack{\\$Z^2,\, W^{28},\,X^7,\, WY^3$}
& $1,4,9,14;28 $
& $1,6,14,21;42$ 
&  \shortstack{\\$Z^2,\, W^{42},\,X^7,\, Y^3$}\\
\hline
$4.$ 
&  \shortstack{\\$Z^2,\, W^{36},\,WX^7,\, Y^3$} 
& $1,5,12,18;36$ 
& $1,6,14,21;42$ 
&  \shortstack{\\$Z^2,\, W^{42},\,X^7,\, Y^3$} \\
\hline
$5.$ 
&  \shortstack{\\$WZ^2,\, W^{21},\,X^7,\, Y^3$} 
& $1,3,7,10;21$ 
& $1,3,7,10;21$ 
&  \shortstack{\\$WZ^2,\, W^{21},\,X^7,\, Y^3$} \\
\hline
$6.$ 
&  \shortstack{\\$Z^2,\, W^{28},\, X^7,\, WY^3$} 
& $1,4,9,14;28$ 
& $1,3,7,10;21$ 
&  \shortstack{\\$WZ^2,\, W^{21},\, X^7,\, Y^3$} \\
\hline
$7.$ 
&  \shortstack{\\$Z^2,\, W^{36},\, WX^7,\, Y^3$} 
& $1,5,12,18;36$ 
& $1,3,7,10;21$ 
&  \shortstack{\\$WZ^2,\, W^{21},\, X^7,\, Y^3$} \\
\hline
$8.$ 
&  \shortstack{\\$Z^2,\, W^{28},\, X^7,\, WY^3$} 
& $1,4,9,14;28$ 
& $1,4,9,14;28$ 
&  \shortstack{\\$Z^2,\, W^{28},\, X^7,\, WY^3$} \\
\hline
$9.$ 
&  \shortstack{\\$Z^2,\, W^{36},\, WX^7,\, Y^3$}
& $1,5,12,18;36$ 
& $1,4,9,14;28$ 
&  \shortstack{\\$Z^2,\, W^{28},\, X^7,\, WY^3$}\\
\hline
$10.$ 
&  \shortstack{\\$Z^2,\, W^{36},\, WX^7,\, Y^3$}
& $1,5,12,18;36$ 
& $1,5,12,18;36$ 
&  \shortstack{\\$Z^2,\, W^{36},\, WX^7,\, Y^3$}\\
\hline
\hline
\raisebox{3mm}{$11.$} 
& \shortstack{\\ $Z^2,\, W^{30},\, W^6X^6,$ $X^5Y,\, Y^3$ \\ \dotfill \\ $Z^2,\, W^{30},\, W^2X^7,$ $X^5Y,\, Y^3$} 
& \raisebox{3mm}{$1,4,10,15;30$}
& \raisebox{3mm}{$1,6,8,15;30$}
& \shortstack{\\ $Z^2,\, W^{30},\, X^5,$ $XY^3,\, W^6Y^3$ \\ \dotfill \\ $Z^2,\, W^{30},\, X^5,$ $XY^3,\, W^{14}Y^2$} \\
\hline
\raisebox{3mm}{$12.$}
& \shortstack{\\ $Z^2,\, W^{30},\, W^6X^6,$ $X^5Y,\, Y^3$ \\ \dotfill \\ $Z^2,\, W^{30},\, W^2X^7,$ $X^5Y,\, Y^3$} 
& \raisebox{3mm}{$1,4,10,15;30$}
& \raisebox{3mm}{$1,5,7,13;26$} 
& \shortstack{\\$Z^2,\, W^{26},\, WX^5,$ $XY^3,\, W^5Y^3$ \\ \dotfill \\ $Z^2,\, W^{26},\, WX^5,$ $XY^3,\, W^{12}Y^2$} \\
\hline
\raisebox{3mm}{$13.$}
& \shortstack{\\$Z^2,\, W^{30},\, X^5,$ $XY^3,\, W^6Y^3$ \\ \dotfill \\ $Z^2,\, W^{30},\, X^5,$ $XY^3,\, W^{14}Y^2$} 
& \raisebox{3mm}{$1,6,8,15;30$ }
& \raisebox{3mm}{$1,3,7,11;22$}
& \shortstack{\\$Z^2,\, W^{22},\, W^4X^6,$ $X^5Y,\, WY^3$ \\ \dotfill \\ $Z^2,\, W^{22},\, WX^7,$ $X^5Y,\, WY^3$ } \\
\hline
\raisebox{3mm}{$14.$}
 & \shortstack{\\$Z^2,\, W^{26},\, WX^5,$ $XY^3,\, W^5Y^3$\\ \dotfill \\ $Z^2,\, W^{26},\, WX^5,$ $XY^3,\, W^{12}Y^2$}
 & \raisebox{3mm}{$1,5,7,13;26$}
 & \raisebox{3mm}{$1,3,7,11;22$}
 & \shortstack{\\$Z^2,\, W^{22},\, W^6X^6,$ $X^5Y,\, WY^3$ \\ \dotfill \\ $Z^2,\, W^{22},\, X^7,$ $X^5Y,\, WY^3$} \\
\hline
\hline
\raisebox{3mm}{$15.$}
& \shortstack{\\$Z^2,\, W^{24},\, W^6X^6,$ $X^4Z,\, Y^3$\\ \dotfill \\ $Z^2,\, W^{24},\, X^8,\, Y^3$} 
& \raisebox{3mm}{$1,3,8,12;24$}
& \raisebox{3mm}{$1,6,8,9;24$}
& \shortstack{\\$W^6Z^2,\, W^{24},\, X^4,$ $XZ^2,\, Y^3$ \\ \dotfill \\ $XZ^2,\, W^{24},\, X^4,,\, Y^3$} 
\\
\hline
\raisebox{3mm}{$16.$}
& \shortstack{\\$Z^2,\, W^{24},\, W^6X^6,$ $X^4Z,\, Y^3$ \\ \dotfill \\ $Z^2,\, W^{24},\, X^8,\, Y^3$}
& \raisebox{3mm}{$1,3,8,12;24$}
& \raisebox{3mm}{$1,5,7,8;21$}
& \shortstack{\\$W^5Z^2,\, W^{21},\, WX^4,$ $XZ^2,\, Y^3$ \\ \dotfill \\ $XZ^2,\, W^{21},\, WX^4,\, Y^3$}
\\
\hline
$17.$
& \shortstack{\\$W^6Z^2,\, W^{24},\, X^4,$ $XZ^2,\, Y^3$}
& $1,6,8,9;24$
& $1,2,5,7;15$
& \shortstack{\\$WZ^2,\, W^{15},\, W^3X^6,$ $X^4Z,\, Y^3$}
\\
\hline
$18.$
& \shortstack{\\$W^5Z^2,\, W^{21},\, WX^4,$ $XZ^2,\, Y^3$}
& $1,5,7,8;21$
& $1,2,5,7;15$
& \shortstack{\\$WZ^2,\, W^{15},\, W^3X^6,$ $X^4Z,\, Y^3$}
\\
\hline
\hline
\raisebox{10mm}{$19.$}
& \shortstack{\\$Z^2,\, W^{22},\, W^2X^5,$\\ $X^4Y,\, XY^3,\, W^{10}Y^2$ \\ \dotfill 
\\ $Z^2,\, W^{22},\, W^6X^4,$\\ $X^4Y,\, XY^3,\, W^4Y^3$ \\ \dotfill 
\\ $Z^2,\, W^{22},\, W^6X^4,$\\ $X^4Y,\, XY^3,\, W^{10}Y^2$} 
& \raisebox{10mm}{$ 1,4,6,11;22$}
& \raisebox{10mm}{$1,4,6,11;22$}
& \shortstack{\\$Z^2,\, W^{22},\, W^2X^5,$\\ $X^4Y,\, XY^3,\, W^{10}Y^2$  \\ \dotfill 
\\ $Z^2,\, W^{22},\, W^6X^4,$\\ $X^4Y,\, XY^3,\, W^4Y^3 $\\ \dotfill 
\\ $Z^2,\, W^{22},\, W^2X^5,$\\ $X^4Y,\, XY^3,\, W^4Y^3$}
\\
\hline
\hline
\raisebox{3mm}{$20.$}
& \shortstack{\\$Z^2,\, W^{18},\, X^6,$ $XY^3,\, W^8Y^2$ \\ \dotfill \\ $Z^2,\, W^{18},\, X^6,$ $XY^3,\, W^2Y^3$}
& \raisebox{3mm}{$1,3,5,9;18$}
& \raisebox{3mm}{$1,4,6,7;18$}
& \shortstack{\\$XZ^2,\, W^{18},\, W^2X^4,$ $X^3Y,\, Y^3$ \\ \dotfill \\ $XZ^2,\, W^{18},\, W^6X^3,$ $X^3Y,\, Y^3$}
\\
\hline
\hline
\raisebox{3mm}{$21.$}
& \shortstack{\\ $XZ^2,\, W^{15},\, X^5,\, Y^3$ \\ \dotfill \\ $W^3Z^2,\, W^{15},\, X^5,$ $XZ^2,\, Y^3$}
& \raisebox{3mm}{$1,3,5,6;15$}
& \raisebox{3mm}{$1,3,5,6;15$}
& \shortstack{\\$XZ^2,\, W^{15},\, X^5,\, Y^3$  \\ \dotfill \\ $XZ^2,\, W^{15},\, W^6X^3,$ $X^3Z,\, Y^3$}
\\
\hline
\hline
$22.$ 
&  \shortstack{\\$Z^2,\, W^{20},\, X^5,\, Y^2Z$}
& $1,4,5,10;20$ 
& $1,4,5,10;20$  
&  \shortstack{\\$Z^2,\, W^{20},\, X^5,\, Y^4$}\\
\hline
$23.$ 
&  \shortstack{\\$Z^2,\, W^{20},\, X^5,\, Y^4$} 
& $1,4,5,10;20$ 
& $1,3,4,7;15$ 
&  \shortstack{\\$WZ^2,\, W^{15},\, X^5,\, Y^2Z$}\\
\hline
$24.$ 
&  \shortstack{\\Not exist} 
& $1,3,4,7;15$ 
& $1,3,4,7;15$ 
&  \shortstack{\\Not exist}
\\
\hline
\hline
\raisebox{3mm}{$25.$}
& \shortstack{\\ $Z^2,\, W^{16},\, WX^5,$ $X^4Y,\, Y^4$ \\ \dotfill \\ $Z^2,\, W^{16},\, W^4X^4,$ $X^4Y,\, Y^4$}
& \raisebox{3mm}{$1,3,4,8;16$}
& \raisebox{3mm}{$1,4,5,6;16$}
& \shortstack{\\ $Y^2Z,\, W^{10}Z,\, W^{16},$ $X^4,\, XZ^2$ \\ \dotfill \\ $WY^3,\, XZ^2,\, W^4Z^2,$ $W^{16},\, X^4$}
\\
\hline
\hline
\raisebox{14mm}{$26.$}
& \shortstack{\\ $W^3Z^2,\, Y^2Z,\, XZ^2,\, W^{13},$\\ $W^4X^3,\, X^3Y,\,WY^3$ \\ \dotfill \\ $W^8Z,\, Y^2Z,\, XZ^2,\,W^{13},$\\ $WX^4,\, X^3Y,\,WY^3$ \\ \dotfill \\ $W^3Z^2,\, Y^2Z,\, XZ^2,\,W^{13},$\\ $WX^4,\, X^3Y,$ $W^9Y$  \\ \dotfill \\ $W^3Z^2,\, Y^2Z,\, XZ^2,\,W^{13},$\\ $WX^4,\, X^3Y,\,WY^3$}
& \raisebox{14mm}{$1,3,4,5;13$}
& \raisebox{14mm}{$1,3,4,5;13$}
& \shortstack{\\ $W^3Z^2,\, Y^2Z,\, XZ^2,\,W^{13},$\\ $W^4X^3,\, X^3Y,$ $W^9Y$ \\ \dotfill \\ $W^8Z,\, Y^2Z,\, XZ^2,\,W^{13},$\\ $WX^4,\, X^3Y,$ $W^9Y$ \\ \dotfill \\  $W^8Z,\, Y^2Z,\, XZ^2,\,W^{13},$\\ $W^4X^3,\, X^3Y,\,WY^3$ \\ \dotfill \\ $W^8Z,\, Y^2Z,\, XZ^2,\,W^{13},$\\ $W^4X^3,\, X^3Y,$ $W^9Y$}
\\
\hline
\hline
$27.$ 
&  \shortstack{\\Not exist} 
& $1,3,4,4;12$ 
& $1,3,4,4;12$  
&  \shortstack{\\Not exist} \\
\hline
\hline
\raisebox{3mm}{$28.$}
& \shortstack{\\ $Z^2,\, W^{18},\, X^9,\, Y^3$ \\ \dotfill \\ $Z^2,\, W^{18},\, W^2X^8,\, Y^3$}
& \raisebox{3mm}{$1,2,6,9;18$}
& \raisebox{3mm}{$2,3,8,11;24$}
& \shortstack{\\$WZ^2,\, W^9X^2,\, X^8,\, Y^3$ \\ \dotfill \\ $WZ^2,\, W^{12},\, X^8,\, Y^3$}
\\
\hline
\raisebox{3mm}{$29.$}
& \shortstack{\\$Z^2,\, W^{18},\, X^9,\, Y^3$ \\ \dotfill \\ $Z^2,\, W^{18},\, W^2X^8,\, Y^3$}
& \raisebox{3mm}{$1,2,6,9;18$}
&\raisebox{3mm}{ $2,5,14,21;42$}
& \shortstack{\\$Z^2,\, W^{16}X^2,\, WX^8,\, Y^3$ \\ \dotfill \\ $Z^2,\, W^{21},\, WX^8,\, Y^3$}
\\
\hline
\hline
\raisebox{3mm}{$30.$}
& \shortstack{\\$Z^2,\, W^{14},\, X^7,\, XY^3,\, W^2Y^3$ \\ \dotfill \\ $Z^2,\, W^{14},\, X^7,\,XY^3,\, W^6Y^2$}
& \raisebox{3mm}{$1,2,4,7;14$}
& \raisebox{3mm}{$2,3,8,13;26$}
& \shortstack{\\$Z^2,\, W^{10}X^2,\, W^4X^6,\,X^6Y,\, WY^3$ \\ \dotfill \\ $Z^2,\, W^{10}X^2,\, WX^8,\,X^6Y,\, WY^3$}
 \\
\hline
\hline
\raisebox{10mm}{$31.$}
& \shortstack{\\$XZ^2,\, W^{12},\, W^2X^5,\, Y^3$ \\ \dotfill \\ $XZ^2,\, W^{12},\, X^6,\, Y^3$ \\ \dotfill \\ $Y^3,\, W^2Z^2,\, W^{12},\,X^6,\, XZ^2$ \\ \dotfill \\ $Y^3,\, XZ^2,\, X^4Y,$\\ $W^2X^5,\, W^{12},\, W^2Z^2$}
& \raisebox{10mm}{$1,2,4,5;12$}
& \raisebox{10mm}{$2,3,10,15;30$}
& \shortstack{\\$Z^2,\, W^{15},\, X^{10},\, Y^3$ \\ \dotfill \\ $Z^2,\, W^{12}X^2,\, X^{10},\, Y^3$ \\ \dotfill \\ $Y^3,\, Z^2,\, W^{12}X^2,\,W^6X^6,\, X^5Z$ \\ \dotfill \\ $Y^3,\, Z^2,\, W^5Y^2,$\\$W^{12}X^2,\, W^6X^6,\, X^5Z$}
\\
\hline
\hline
\raisebox{8mm}{$32.$}
& \shortstack{\\$Z^2,\, W^{10},\, X^5,\, Y^5$ \\ \dotfill \\ $Z^2,\, W^{10},\, X^5,$\\  $W^2Y^4,\, X^{a+1}Y^{-a+4}$\\ $a=0,1,2,3. $}
& \raisebox{8mm}{$1,2,2,5;10$}
& \raisebox{8mm}{$2,4,5,9;20$}
& \shortstack{\\$WZ^2,\, W^5Y^2,\, X^5,\, Y^4$ \\ \dotfill \\ $WZ^2,\, W^5Y^2,\, X^5,$\\  $Y^4,\, W^{2(a+1)}X^{-a+4}$\\ $a=0,1,2,3. $}
\\
\hline
\raisebox{8mm}{$33.$}
& \shortstack{\\ $Z^2,\, W^{10},\, X^5,\, Y^5$ \\ \dotfill \\ $Z^2,\, W^{10},\, X^5,$\\ $W^2Y^4,\, X^{a+1}Y^{-a+4}$ \\ $a=0,1,2,3. $}
& \raisebox{8mm}{$1,2,2,5;10$}
& \raisebox{8mm}{$2,6,7,15;30$}
& \shortstack{\\$Z^2,\, W^5Y^2,\, X^5,\, WY^4$ \\ \dotfill \\ $Z^2,\, W^8Y^2,\, X^5,$\\ $WY^4,\, W^{3(a+1)}X^{-a+4}$\\ $a=0,1,2,3. $}
\\
\hline
\raisebox{8mm}{$34.$}
& \shortstack{\\$Z^2,\, W^{10},\, X^5,\, Y^5$ \\ \dotfill \\ $Z^2,\, W^{10},\, X^5,$\\ $W^2Y^4,\, X^{a+1}Y^{-a+4}$ \\ $a=0,1,2,3. $}
& \raisebox{8mm}{$1,2,2,5;10$}
& \raisebox{8mm}{$2,5,6,13;26$}
& \shortstack{\\$Z^2,\, W^8X^2,\, X^4Y,\, WY^4$ \\ \dotfill \\ $Z^2,\, W^8X^2,\, X^4Y,$\\  $WY^4,\, W^{3a+4}Y^{-a+3}$ \\ $a=0,1,2,3. $}
\\
\hline
\hline
$35.$ 
&  \shortstack{\\$Z^2,\, W^{12},\, X^{12},\, Y^3$} 
& $1,1,4,6;12$ 
& $3,5,11,14;33$ 
&  \shortstack{\\$XZ^2,\, W^{11},\, WX^6,\, Y^3$}\\
\hline
$36.$ 
&  \shortstack{\\$Z^2,\, W^{12},\, X^{12},\, Y^3$} 
& $1,1,4,6;12$ 
& $4,5,13,22;44$ 
&  \shortstack{\\$Z^2,\, W^{11},\, WX^8,\, XY^3$}\\
\hline
$37.$ 
&  \shortstack{\\$Z^2,\, W^{12},\, X^{12},\, Y^3$} 
& $1,1,4,6;12$ 
& $5,6,22,33;66$ 
&  \shortstack{\\$Z^2,\, W^{12}X,\, X^{11},\, Y^3$}\\
\hline
\hline
\raisebox{3mm}{$38.$}
& \shortstack{\\$Z^2,\, W^{10},\, X^{10},\,XY^3,\, WY^3$ \\ \dotfill \\ $Z^2,\, W^{10},\, X^{10},\,XY^3,\, W^4Y^2$}
& \raisebox{3mm}{$1,1,3,5;10$}
& \raisebox{3mm}{$3,4,10,13;30$}
& \shortstack{\\$XZ^2,\, W^{10},\, W^6X^3,\,X^5Y,\, Y^3$ \\ \dotfill \\ $XZ^2,\, W^{10},\, W^2X^6,\,X^5Y,\, Y^3 $}
\\
\hline
$39.$ 
&  \shortstack{\\Not exist} 
& $1,1,3,5;10$ 
& $3,5,11,19;38$ 
&  \shortstack{\\Not exist}\\
\hline
\raisebox{3mm}{$40.$}
& \shortstack{\\$Z^2,\, W^{10},\, X^{10},\,XY^3,\, WY^3$ \\ \dotfill \\ $Z^2,\, W^{10},\, X^{10},\,X^4Y^2,\, WY^3$}
& \raisebox{3mm}{$1,1,3,5;10$}
& \raisebox{3mm}{$4,5,18,27;54$}
& \shortstack{\\$Z^2,\, W^9Y,\, W^6X^6,\,WX^{10},\, Y^3$ \\ \dotfill \\ $Z^2,\, W^{11}X^2,\, WX^{10},\,Y^3,\, W^9Y$}
\\
\hline
\hline
\raisebox{3mm}{$41.$}
& \shortstack{\\$WZ^2,\, W^9,\, X^9,\, Y^3$ \\ \dotfill \\ Y$^3,\, WZ^2,\, W^9,\,X^9,\, XZ^2$}
& \raisebox{3mm}{$1,1,3,4;9$}
& \raisebox{3mm}{$3,4,11,18;36$}
& \shortstack{\\$Z^2,\, W^{12},\, X^9,\, WY^3$ \\ \dotfill \\ $WY^3,\, Z^2,\, W^6Z,\,W^4X^6,\, X^9$}
\\
\hline
$42.$ 
&  \shortstack{\\$Y^3,\, WZ^2,\, W^9,\,X^9,\, XZ^2$} 
& $1,1,3,4;9$ 
& $2,5,9,11;27$ 
&  \shortstack{\\$Y^3,\, XZ^2,\, W^8Z,\,W^6X^3,\, WX^5$}\\
\hline
\raisebox{3mm}{$43.$}
& \shortstack{\\$WZ^2,\, W^9,\, X^9,\, Y^3$ \\ \dotfill \\ $Y^3,\, WZ^2,\, W^9,\,X^9,\, XZ^2$}
& \raisebox{3mm}{$1,1,3,4;9$}
& \raisebox{3mm}{$3,5,16,24;48$}
& \shortstack{\\$Z^2,\, W^{16},\, WX^9,\, Y^3$ \\ \dotfill \\ $Y^3,\, Z^2,\, W^8Z,\,W^6X^6,\, WX^9$}
\\
\hline
\hline
$44.$ 
&  \shortstack{\\$Z^2,\, W^8,\, X^8,\, Y^4$}
& $1,1,2,4;8$ 
& $3,4,7,10;24$ 
&  \shortstack{\\$XZ^2,\, W^8,\, X^6,\, Y^2Z$}\\
\hline
\hline
$45.$
&  \shortstack{\\$Y^2Z,\, XZ^2,\, WZ^2,\,W^7,\, X^7$}
& $1,1,2,3;7$
& $3,4,7,14;28$
&  \shortstack{\\$Z^2,\, W^7Y,\, W^4X^4,\,X^7,\, Y^4$}
\\
\hline
\hline
\raisebox{3mm}{$46.$}
& \shortstack{\\$Y^3Z,\, WZ^2,\, X^3Z,\, W^5,\, X^5,$\\ $Y^5$}
& \raisebox{3mm}{$1,1,1,2;5$}
& \raisebox{3mm}{$4,5,7,9;25$}
& \shortstack{\\$YZ^2,\, W^4Z,\, W^5X,\, WY^3,$\\ $X^5$}
\\
\hline
$47.$ 
&  \shortstack{\\$WZ^2,\, W^5,\, X^5,\, Y^5$} 
& $1,1,1,2;5$  
& $5,7,8,20;40$ 
&  \shortstack{\\$Z^2,\, W^8,\, WX^5,\, Y^5$}\\
\hline
\hline
$48.$ 
&  \shortstack{\\$Z^2,\, W^6,\, X^6,\, Y^6$} 
& $1,1,1,3;6$ 
& $5,6,8,11;30$  
&  \shortstack{\\$YZ^2,\, W^6,\, X^5,\, XY^3$}\\
\hline
$49.$ 
&  \shortstack{\\$Z^2,\, W^6,\, X^6,\, Y^6$} 
& $1,1,1,3;6$ 
& $7,8,10,25;50$ 
&  \shortstack{\\$Z^2,\, W^6X,\, X^5Y,\, Y^5$}\\
\hline
\hline
$50.$ 
&  \shortstack{\\$Z^4,\, W^4,\, X^4,\, Y^4$} 
& $1,1,1,1;4$ 
& $7,8,9,12;36$  
&  \shortstack{\\$Z^3,\, W^4X,\, X^3Z,\, Y^4$}\\
\hline
\hline
$51.$
&  \shortstack{\\$Z^2,\, W^7,\, X^7,\,XY^4,\, WY^4$}
& $2,2,3,7;14$
& $2,2,3,7;14$
&  \shortstack{\\$Z^2,\, W^4Y^2,\, W^7,\,X^4Y^2,\, XY^4$}
\\
\hline
\caption{Polytope duality associated to coupling pairs}\label{ListMainThm}
\end{longtable}
\endgroup
\end{thm}

\begin{remark}
In Table \ref{ListMainThm}, the reflexive polytopes $\Delta$ and $\Delta'$ are given as a set of monomials that are vertices of them. 
If there are more than one pairs, they are separated by a dotted line and polytopes in the same row give the polytope duality. 
\end{remark}

\bigskip
\noindent
\proof
Take polynomials $F$ and $F'$ that are respectively anticanonical sections of the weighted projective spaces $\mathbb{P}(a)$ and $\mathbb{P}(b)$ as in Table ~\ref{ListMainThm} in each case. 

Recall that a pair of reflexive polytopes $\Delta$ and $\Delta'$ is {\it polytope dual} if relations $\Delta_F\subset\Delta\subset\Delta_a,\, \Delta_{F'}\subset\Delta'\subset\Delta_b$, and $\Delta^*\simeq\Delta'$ hold. 

The strategy of the proof is that in each case, after taking a basis of the lattice $M_n$,  we observe if the Newton polytope $\Delta_F$ of the polynomial $F$ is reflexive by a direct computation. 
If the polytope is not reflexive, then, we search a reflexive polytope $\Delta$ satisfying inclusions $\Delta_F\subset\Delta\subset\Delta^{(n)}$ of polytopes. 
The analogous observations should be made for $F'$. 
Once one gets a candidate reflexive polytope $\Delta$ and $\Delta'$, we then study whether they satisfy a relation $\Delta^*\simeq\Delta'$. 

The assertion is proved case by case. 

\subsection{No. 1--No. 10} 
We claim that the unique pair $(\Delta^{(14)},\, {\Delta^{(14)}}^*)$ is polytope dual commonly for Nos. 1 to 10. 
Take a basis $\{ e_1^{(n)},\, e_2^{(n)},\, e_3^{(n)}\}$ of a lattice $M_{n}$ for $n=14,28,45,$ and $51$ by 
\begin{eqnarray*}
e_1^{(14)} = (-6,1,0,0), & &e_2^{(14)} = (-14,0,1,0),\quad e_3^{(14)} = (-21,0,0,1),\, \\
e_1^{(28)} = (-3,1,0,0), & &e_2^{(28)} = (-7,0,1,0),\quad e_3^{(28)} = (-10,0,0,1),\, \\
e_1^{(45)} = (-4,1,0,0), & &e_2^{(45)} = (-9,0,1,0),\quad e_3^{(45)} = (-14,0,0,1),\, \\ 
e_1^{(51)} = (-5,1,0,0),& & e_2^{(51)} = (-12,0,1,0),\quad e_3^{(51)} = (-18,0,0,1).   
\end{eqnarray*}

In \cite{KM}, it is proved that the polytopes $\Delta^{(n)}$ for $n=14,28,45,$ and $51$ are isomorphic to the polytope $\Delta^{(14)}$ that is the convex hull of vertices $(-1,-1,1)$,\, $(-1,-1,-1)$,\, $(6,-1,-1)$, and $(-1,2,-1)$ under the above choice of basis. 
Since the polar dual ${\Delta^{(14)}}^*$ is the convex hull of vertices $(1,0,0)$,\, $(0,1,0)$,\, $(0,0,1)$, and $(-6,-14,-21)$, the linear map of $\mathbb{R}^3$ defined by a matrix $\left(\begin{smallmatrix} 6 & -1 & -1\\ -1 & 2 & -1\\ -1 & -1 & 1 \end{smallmatrix}\right)$ gives an isomorphism from $\Delta^{(14)}$ to ${\Delta^{(14)}}^*$. 
Therefore, a relation ${\Delta^{(14)}}^*\simeq\Delta^{(14)}$ holds.  

Define $2$-dimensional faces $\Gamma_1,\, \Gamma_2$, and $\Gamma_3$ by 
\begin{eqnarray*}
\Gamma_1 & = & {\rm Conv}\{ (-1,-1,0),\, (6,-1,-1),\, (-1,2,-1)\}, \\
\Gamma_2 & = & {\rm Conv}\{ (-1,0,-1),\, (6,-1,-1),\, (-1,-1,1)\}, \\
\Gamma_3 & = & {\rm Conv}\{ (0,-1,-1),\, (-1,-1,1),\, (-1,2,-1)\}. 
\end{eqnarray*}

\begin{lem}\label{LemNonRef}
If a subpolytope $\Delta$ of $\Delta^{(14)}$ contains at least one of the faces $\Gamma_1,\, \Gamma_2$, and $\Gamma_3$, then, $\Delta$ is not reflexive, and there does not exist a reflexive polytope $\tilde{\Delta}$ such that $\Delta\subset\tilde{\Delta}\subset\Delta^{(14)}$ except $\Delta^{(14)}$. 
\end{lem}
\proof
The polar duals of the faces $\Gamma_1,\, \Gamma_2$, and $\Gamma_3$ are respectively rational vertices 
$(3/10,\, 7/10,\, 21/10), $\,
$(2/9,\, 14/9,\, 7/9), $ and 
$(6/5,\, 2/5,\, 3/5)$. 
Therefore, $\Delta$ is not reflexive. 
The last assertion follows immediately. 
\QED

It is easy to see that any reflexive subpolytope $\tilde{\Delta}$ of $\Delta^{(14)}$ should contain faces of the form 
\[
\begin{array}{c}
{\rm Conv}\{ (-1,-1,0),\, (6,-1,-1),\, (-1,a,-1)\}, \quad\textnormal{or}\\
{\rm Conv}\{ (-1,-1,1),\, (-1,0,-1),\, (b,-1,-1)\}, \quad\textnormal{or}\\
{\rm Conv}\{ (-1,-1,1),\, (0,-1,-1),\, (-1,c,-1)\}
\end{array}
\]
with $a=-1,2,1$, $b=-1,0,1,2,3,4,5,6$, $c=-1,0,1$, of which the polar duals are respectively vertices of the form 
\[
\frac{1}{a+8}(a+1, \, 7, \, 7(a+1)), \, 
\frac{1}{b+3}(2, \, 2(b+1), \, b+1), \,
\frac{1}{c+3}(2(c+1), \, 2, \, c+1). 
\]
In order that $\tilde{\Delta}$ should be reflexive, one has $a=b=c=-1$, that is, $\tilde{\Delta}$ should be $\Delta^{(14)}$. 

\noindent
{\bf No. 1.}
The Newton polytope of $F=F'$ coincides with $\Delta^{(14)}$. 

\noindent
{\bf No. 2.}
The Newton polytope of $F$, which is the convex hull of vertices $(-1,-1,0)$,\, $(-1,-1,1)$,\, $(6,-1,-1)$, and $(-1,2,-1)$, contains the face $\Gamma_1$, thus by Lemma \ref{LemNonRef}, it is not reflexive, and $\Delta^{(14)}$ is the only reflexive polytope. 
Besides, the Newton polytope of $F'$ of the coupling dual partner , which is the convex hull of vertices $(-1,-1,1)$,\, $(-1,-1,-1)$,\, $(6,-1,-1)$, and $(-1,2,-1)$, coincides with $\Delta^{(28)}$. 

\noindent
{\bf No. 3.}
The Newton polytope of $F$, which is the convex hull of vertices $(-1,0,-1)$,\, $(-1,-1,1)$,\, $(6,-1,-1)$, and $(-1,2,-1)$, contains the face $\Gamma_2$, thus by Lemma \ref{LemNonRef}, it is not reflexive, and $\Delta^{(14)}$ is the only reflexive polytope. 
Besides, the Newton polytope of $F'$, which is the convex hull of vertices $(-1,-1,1)$,\, $(-1,-1,-1)$,\, $(6,-1,-1)$, and $(-1,2,-1)$, coincides with $\Delta^{(45)}$. 

\noindent
{\bf No. 4.}
The Newton polytope of $F$, which is the convex hull of vertices $(0,-1,-1)$,\, $(-1,-1,1)$,\, $(6,-1,-1)$, and $(-1,2,-1)$, contains the face $\Gamma_3$, thus by Lemma \ref{LemNonRef}, it is not reflexive, and $\Delta^{(14)}$ is the only reflexive polytope. 
Besides, the Newton polytope of $F'$, which is the convex hull of vertices $(-1,-1,1)$,\, $(-1,-1,-1)$,\, $(6,-1,-1)$, and $(-1,2,-1)$, coincides with $\Delta^{(51)}$. 

\noindent
{\bf No. 5.}
The Newton polytope of $F=F'$, which is the convex hull of vertices $(-1,-1,0)$,\, $(-1,-1,1)$,\, $(6,-1,-1)$, and $(-1,2,-1)$, contains the face $\Gamma_1$, thus by Lemma \ref{LemNonRef}, it is not reflexive, and $\Delta^{(28)}$ is the only reflexive polytope. 

\noindent
{\bf No. 6.}
The Newton polytope of $F$, which is the convex hull of vertices $(-1,0,-1)$,\, $(-1,-1,1)$,\, $(6,-1,-1)$, and $(-1,2,-1)$, contains the face $\Gamma_2$, thus by Lemma \ref{LemNonRef}, it is not reflexive, and $\Delta^{(28)}$ is the only reflexive polytope. 
Besides, the Newton polytope of $F'$, which is the convex hull of vertices $(-1,-1,0)$,\, $(-1,-1,1)$,\, $(6,-1,-1)$, and $(-1,2,-1)$, contains the face $\Gamma_1$, thus by Lemma \ref{LemNonRef}, it is not reflexive, and $\Delta^{(45)}$ is the only reflexive polytope. 

\noindent
{\bf No. 7.}
The Newton polytope of $F$, which is the convex hull of vertices $(0,-1,-1)$,\, $(-1,-1,1)$,\, $(6,-1,-1)$, and $(-1,2,-1)$, contains the face $\Gamma_3$, by Lemma \ref{LemNonRef}, it is not reflexive, and $\Delta^{(28)}$ is the only reflexive polytope. 
Besides, the Newton polytope of $F'$, which is the convex hull of vertices $(-1,-1,0)$,\, $(-1,-1,1)$,\, $(6,-1,-1)$, and $(-1,2,-1)$, contains the face $\Gamma_1$, thus by Lemma \ref{LemNonRef}, it is not reflexive, and $\Delta^{(51)}$ is the only reflexive polytope. 

\noindent
{\bf No. 8.}
The Newton polytope of $F=F'$, which is the convex hull of vertices $(-1,0,-1)$,\, $(-1,-1,1)$,\, $(6,-1,-1)$, and $(-1,2,-1)$, contain the face $\Gamma_1$, thus by Lemma \ref{LemNonRef}, it is not reflexive, and $\Delta^{(45)}$ is the only reflexive polytope. 

\noindent
{\bf No. 9.}
The Newton polytope of $F$, which is the convex hull of vertices $(0,-1,-1)$,\, $(-1,-1,1)$,\, $(6,-1,-1)$, and $(-1,2,-1)$, contains the face $\Gamma_3$, thus by Lemma \ref{LemNonRef}, it is not reflexive, and $\Delta^{(45)}$ is the only reflexive polytope. 
Besides, the Newton polytope of $F'$, which is the convex hull of vertices $(-1,0,-1)$,\, $(-1,-1,1)$,\, $(6,-1,-1)$, and $(-1,2,-1)$, contains the face $\Gamma_2$, thus by Lemma \ref{LemNonRef}, it is not reflexive, and $\Delta^{(51)}$ is the only reflexive polytope. 

\noindent
{\bf No. 10.}
The Newton polytope of $F=F'$, which is the convex hull of vertices $(0,-1,-1)$,\, $(-1,-1,1)$,\, $(6,-1,-1)$, and $(-1,2,-1)$, contains the face $\Gamma_3$, thus by Lemma \ref{LemNonRef}, they are not reflexive, and $\Delta^{(51)}$ is the only reflexive polytope. 

\subsection{No. 11--No. 14} 
We claim that there exist two polytope-dual pairs for Nos. 11 to 14. 
Take a basis $\{ e_1^{(n)},\, e_2^{(n)},\, e_3^{(n)}\}$ for a lattice $M_{n}$ for $n=38,50,77,$ and $82$ by 
\begin{eqnarray*}
e_1^{(38)} = (-6,1,0,0),& & e_2^{(38)} = (-8,0,1,0),\quad e_3^{(38)} = (-15,0,0,1),\, \\
e_1^{(50)} = (-4,1,0,0),& & e_2^{(50)} = (-10,0,1,0),\quad e_3^{(50)} = (-15,0,0,1),\, \\
e_1^{(77)} = (-5,1,0,0),& & e_2^{(77)} = (-7,0,1,0),\quad e_3^{(77)} = (-13,0,0,1),\, \\ 
e_1^{(82)} = (-3,1,0,0),& & e_2^{(82)} = (-7,0,1,0),\quad e_3^{(82)} = (-11,0,0,1).   
\end{eqnarray*}

In \cite{KM}, it is proved that the polytopes $\Delta^{(38)}$ and $\Delta^{(77)}$, and $\Delta^{(50)}$, and $\Delta^{(82)}$ are respectively isomorphic to the polytopes 
\[
\begin{array}{c}
\Delta^{(38,77)}:={\rm Conv}\left\{\begin{array}{l} (-1,-1,1),\, (-1,-1,-1),\\ (4,-1,-1),\, (0,2,-1),\, (-1,2,-1)\end{array}\right\}, \\
\Delta^{(50,82)}:={\rm Conv}\left\{\begin{array}{l} (-1,-1,1),\, (-1,-1,-1),\\ (6,-1,-1),\, (4,0,-1),\, (-1,2,-1)\end{array}\right\}. 
\end{array}
\]

Define polytopes $\Delta_1,\, \Delta_1'$, and $\Delta_2,\, \Delta_2'$ by 
\begin{eqnarray*}
\Delta_1 & := & \Delta^{(38,77)}, \\
\Delta_1' & := & {\rm Conv}\left\{\begin{array}{l} (-1,-1,1),\, (-1,-1,-1),\\ (5,-1,-1),\, (4,0,-1),\, (-1,2,-1)\end{array}\right\}, \\
\Delta_2 & := & {\rm Conv}\left\{\begin{array}{l} (-1,-1,1),\, (-1,-1,-1),\\ (6,-1,-1),\, (4,0,-1),\, (-1,2,-1)\end{array}\right\}, \\
\Delta_2' & := & \Delta^{(50,82)}. 
\end{eqnarray*}
Since the polar dual polytopes $\Delta_1^*$ and $\Delta_2^*$ of $\Delta_1$ and $\Delta_2$ are the convex hulls of vertices $(1,0,0)$,\, $(0,1,0)$,\, $(0,0,1)$,\, $(0,-2,-3)$, and $(-6,-8,-15)$, respectively, $(1,0,0)$,\, $(0,1,0)$,\, $(0,0,1)$,\, $(-2,-4,-7)$, and $(-4,-10,-15)$, and the linear maps of $\mathbb{R}^3$ determined by matrices $A_1:=\left(\begin{smallmatrix} 1 & 1 & 2\\ 2 & 3 & 5\\ 3 & 4 & 8 \end{smallmatrix}\right)$, and $A_2:={}^tA_1$ respectively give isomorphisms from $\Delta_1'$ to $\Delta_1^*$ and from  $\Delta_2'$ to $\Delta_2^*$, the relations $\Delta_1^*\simeq\Delta_1'$ and $\Delta_2^*\simeq\Delta_2'$ hold. 

\noindent
{\bf No. 11.}
The Newton polytope of $F$, which is the convex hull of vertices $(-1,-1,1)$,\, $(-1,-1,-1)$,\, $(4,-1,-1)$, and $(0,2,-1)$, is not reflexive. 
It is observed that one might take a polytope with a vertex $(-1,a,-1)$ with $a$ being $1$ or $2$ instead of a face spanned by vertices $(-1,-1,1)$,\, $(-1,-1,-1)$, and $(0,2,-1)$.  
Besides, the Newton polytope of $F'$, which is the convex hull of vertices $(-1,-1,1)$,\, $(-1,-1,-1)$,\, $(4,0,-1)$, and $(-1,2,-1)$, is not reflexive. 
It is observed that one might take a polytope with a vertex $(b,-1,-1)$ with $b$ being $5$ or $6$ instead of a face spanned by vertices $(-1,-1,1)$,\, $(-1,-1,-1)$, and $(4,0,-1)$.  
Therefore, there are two polytope-dual pairs, that is, if $(a,b)=(2,5)$, then, $(\Delta_1, \, \Delta_1')$, and if $(a,b)=(1,6)$, then, $(\Delta_2',\, \Delta_2)$. 

\noindent
{\bf No. 12.}
The Newton polytope of $F$, which is the convex hull of vertices $(-1,-1,1)$,\, $(-1,-1,-1)$,\, $(4,-1,-1)$, and $(0,2,-1)$, is not reflexive. 
It is observed that one might take a polytope with a vertex $(-1,a,-1)$ with $a$ being $1$ or $2$ instead of a face spanned by vertices $(-1,-1,1)$,\, $(-1,-1,-1)$, and $(0,2,-1)$.  
Besides, the Newton polytope of $F'$, which is the convex hull of vertices $(-1,-1,1)$,\, $(0,-1,-1)$,\, $(4,0,-1)$, and $(-1,2,-1)$, is not reflexive. 
It is observed that one might take a polytope with vertices $(-1,-1,-1)$, and $(b,-1,-1)$ with $b$ being $5$ or $6$, instead of faces spanned by vertices $(-1,-1,1)$,\, $(0,-1,-1)$, \, $(4,0,-1)$, and $(-1,-1,1)$,\, $(0,-1,-1)$, \, $(-1,2,-1)$.  
Therefore, there are two polytope-dual pairs, that is, if $(a,b)=(2,5)$, then, $(\Delta_1, \, \Delta_1')$, and if $(a,b)=(1,6)$, then, $(\Delta_2',\, \Delta_2)$. 

\noindent
{\bf No. 13.}
The Newton polytope of $F$, which is the convex hull of vertices $(-1,-1,1)$,\, $(-1,-1,-1)$,\, $(4,0,-1)$, and $(-1,2,-1)$, is not reflexive. 
It is observed that one might take a polytope with a vertex $(a,-1,-1)$ with $a$ being $5$ or $6$, instead of a face spanned by vertices $(-1,-1,1)$,\, $(-1,-1,-1)$, and $(4,0,-1)$. 
Besides, the Newton polytope of $F'$, which is the convex hull of vertices $(-1,-1,1)$,\, $(-1,0,-1)$,\, $(4,-1,-1)$, and $(0,2,-1)$, is not reflexive. 
It is observed that one might take a polytope with a vertex $(-1,-1,-1)$, and $(-1,b,-1)$ with $b$ being $1$ or $2$, instead of faces spanned by vertices $(-1,-1,1)$,\, $(-1,0,-1)$, \, $(4,-1,-1)$, and $(-1,-1,1)$,\, $(-1,0,-1)$,\, $(0,2,-1)$. 
Therefore, there are two polytope-dual pairs, that is, if $(a,b)=(5,2)$, then, $(\Delta_1', \, \Delta_1)$, and if $(a,b)=(6,1)$, then, $(\Delta_2,\, \Delta_2')$. 

\noindent
{\bf No. 14.}
The Newton polytope of $F$, which is the convex hull of vertices $(0,-1,-1)$,\, $(4,0,-1)$,\, $(-1,2,-1)$, and $(-1,-1,1)$, is not reflexive. 
It is observed that one might take a polytope with vertices $(-1,-1,-1)$, and $(a,-1,-1)$ with $a$ being $5$ or $6$, instead of faces spanned by vertices $(-1,-1,1)$,\, $(0,-1,-1)$, \, $(4,0,-1)$, and $(-1,-1,1)$,\, $(0,-1,-1)$,\, $(-1,2,-1)$. 
Besides, the Newton polytope of $F'$, which is the convex hull of vertices $(-1,-1,1)$,\, $(-1,0,-1)$,\, $(4,-1,-1)$, and $(0,2,-1)$, is not reflexive. 
It is observed that one might take a polytope with vertices $(-1,-1,-1)$, and $(-1,b,-1)$ with $b$ being $1$ or $2$, instead of faces spanned by vertices $(-1,-1,1)$,\, $(-1,0,-1)$, \, $(4,-1,-1)$, and $(-1,-1,1)$,\, $(-1,0,-1)$,\, $(0,2,-1)$. 
Therefore, there are two polytope-dual pairs, that is, if $(a,b)=(5,2)$, then, $(\Delta_1', \, \Delta_1)$, and if $(a,b)=(6,1)$, then, $(\Delta_2,\, \Delta_2')$. 

\subsection{No. 15--No. 18} 
We claim that there exists a unique polytope-dual pair for Nos. 17, 18, and that two pairs for Nos. 15, 16. 
Take a basis $\{ e_1^{(n)},\, e_2^{(n)},\, e_3^{(n)}\}$ of a lattice $M_{n}$ for $n=13,20,59,$ and $72$ by 
\begin{eqnarray*}
e_1^{(13)} = (-3,1,0,0), & & e_2^{(13)} = (-8,0,1,0),\quad e_3^{(13)} = (-12,0,0,1),\, \\
e_1^{(20)} = (-6,1,0,0),& & e_2^{(20)} = (-8,0,1,0),\quad e_3^{(20)} = (-9,0,0,1),\, \\
e_1^{(59)} = (-5,1,0,0), & & e_2^{(59)} = (-7,0,1,0),\quad e_3^{(59)} = (-8,0,0,1),\, \\ 
e_1^{(72)} = (-2,1,0,0),& & e_2^{(72)} = (-5,0,1,0),\quad e_3^{(72)} = (-7,0,0,1).   
\end{eqnarray*}

In \cite{KM}, it is proved that the polytopes $\Delta^{(20)}$ and $\Delta^{(59)}$ are isomorphic to the convex hull $\Delta^{(20,59)}$ of vertices $(-1,-1,1)$,\, $(-1,-1,-1)$,\, $(3,-1,-1)$,\, $(0,-1,1)$, and $(-1,2,-1)$. 

Define polytopes $\Delta_1,\, \Delta_1'$, and $\Delta_2,\, \Delta_2'$ by 
\begin{eqnarray*}
\Delta_1 & := & \Delta^{(20,59)}, \\
\Delta_1' & := & {\rm Conv}\left\{\begin{array}{l} (-1,-1,1),\, (-1,-1,-1),\\ (5,-1,-1),\, (3,-1,0),\, (-1,2,-1)\end{array}\right\}, \\
\Delta_2 & := & {\rm Conv}\{ (0,-1,1),\, (-1,-1,-1),\, (3,-1,-1),\, (-1,2,-1)\}, \\
\Delta_2' & := & \Delta^{(13)}={\rm Conv}\{ (-1,-1,1),\, (-1,-1,-1),\, (7,-1,-1),\, (-1,2,-1)\}. 
\end{eqnarray*}
Since the polar dual polytopes $\Delta_1^*$ and $\Delta_2^*$ of $\Delta_1$ and $\Delta_2$ are the convex hulls of vertices $(1,0,0)$,\, $(0,1,0)$,\, $(0,0,1)$,\, $(0,-2,-3)$, and $(-6,-8,-9)$, respectively, $(0,1,0)$,\, $(0,0,1)$,\, $(2,0,-1)$, and $(-6,-8,-9)$, and the linear maps of $\mathbb{R}^3$ determined by a matrix $\left(\begin{smallmatrix} 1 & 1 & 1\\ 2 & 3 & 3\\ 3 & 4 & 5 \end{smallmatrix}\right)$ 
gives isomorphisms from $\Delta_1'$ to $\Delta_1^*$ and from $\Delta_2'$ to $\Delta_2^*$, the relations $\Delta_1^*\simeq\Delta_1'$ and $\Delta_2^*\simeq\Delta_2'$ hold. 

\noindent
{\bf No. 15.}
The Newton polytope of $F$, which is the convex hull of vertices $(0,-1,1)$,\, $(-1,-1,-1)$,\, $(3,-1,-1)$, and $(-1,2,-1)$, coincides with $\Delta_2$. 
Besides, the Newton polytope of $F'$, which is the convex hull of vertices $(-1,-1,1)$,\, $(-1,-1,-1)$,\, $(3,-1,0)$, and $(-1,2,-1)$, is not reflexive. 
It is observed that one might take a polytope with a vertex $(b,-1,-1)$ with $b$ being $5$ or $7$ instead of a face spanned by vertices $(3,-1,0)$,\, $(-1,-1,-1)$, and $(-1,2,-1)$. 
Therefore, there are two polytope-dual pairs, that is, if $b=5$, then, $(\Delta^{(20)} = \Delta_1,\, \Delta_1')$, and if $b=7$, then, $(\Delta_2,\, \Delta^{(13)}=\Delta_2')$. 

\noindent
{\bf No. 16.}
The Newton polytope of $F$, which is the convex hull of vertices $(0,-1,1)$,\, $(-1,-1,-1)$,\, $(3,-1,-1)$, and $(-1,2,-1)$, coincides with $\Delta_2$. 
Besides, the Newton polytope of $F'$, which is the convex hull of vertices $(-1,-1,1)$,\, $(0,-1,-1)$,\, $(3,-1,0)$, and $(-1,2,-1)$, is not reflexive. 
It is observed that one might take a polytope with vertices $(-1,-1,-1)$, and $(-1,b,-1)$ with $b$ being $5$ or $7$ instead of faces spanned by vertices $(-1,-1,1)$,\, $(0,-1,-1)$,\, $(-1,2,-1)$, and $(3,-1,0)$,\, $(0,-1,-1)$,\, $(-1,2,-1)$. 
Therefore, there are two polytope-dual pairs, that is, if $b=5$, then, $(\Delta_1,\, \Delta_1')$, and if $b=7$, then, $(\Delta_2,\, \Delta_2')$. 

\noindent
{\bf No. 17.}
The Newton polytope of $F$, which is the convex hull of vertices $(-1,-1,1)$,\, $(-1,-1,-1)$,\, $(3,-1,0)$, and $(-1,2,-1)$, is not reflexive. 
It is observed that one might take $\Delta^{(72)}$ or a polytope with a vertex $(5,-1,-1)$  instead of a face spanned by vertices $(3,-1,0)$,\, $(-1,-1,-1)$, and $(-1,2,-1)$.  
Besides, the Newton polytope of $F'$, which is the convex hull of vertices $(0,-1,1)$,\, $(-1,-1,0)$,\, $(3,-1,-1)$, and $(-1,2,-1)$, is not reflexive. 
It is observed that one might take a polytope with vertices $(-1,-1,1)$, and $(-1,-1,-1)$ instead of faces spanned by vertices $(0,-1,1)$,\, $(-1,-1,0)$,\, $(-1,2,-1)$, and $(-1,-1,0)$,\, $(3,-1,-1)$,\, $(-1,2,-1)$. 
Therefore, the pair $(\Delta_1',\, \Delta_1)$ is polytope-dual. 

\noindent
{\bf No. 18.}
The Newton polytope of $F$, which is the convex hull of vertices $(-1,-1,1)$,\, $(0,-1,-1)$,\, $(3,-1,0)$, and $(-1,2,-1)$, is not reflexive. 
It is observed that one might take $\Delta^{(72)}$ or a polytope with vertices $(-1,-1,-1)$ and $(5,-1,-1)$ instead of faces spanned by vertices $(0,-1,-1)$,\, $(-1,-1,1)$,\, $(-1,2,-1)$, and $(0,-1,-1)$,\, $(3,-1,0)$,\, $(-1,2,-1)$. 
Besides, the Newton polytope of $F'$, which is the convex hull of vertices $(-1,-1,0)$,\, $(0,-1,1)$,\, $(3,-1,-1)$, and $(-1,2,-1)$, is not reflexive. 
It is observed that one might take a polytope with vertices $(-1,-1,1)$, and $(-1,-1,-1)$ instead of faces spanned by vertices $(-1,-1,0)$,\, $(0,-1,1)$,\, $(-1,2,-1)$, and $(-1,-1,0)$,\, $(-1,2,-1)$,\, $(3,-1,-1)$. 
Therefore, the pair $(\Delta_1',\, \Delta_1)$ is polytope-dual. 

\subsection{No. 19} 
We claim that there exist three polytope-dual pairs. 
Take a basis $\{ e_1^{(78)},\, e_2^{(78)},\, e_3^{(78)}\}$ for a lattice $M_{78}$ by 
\begin{eqnarray*}
e_1^{(78)} = (-4,1,0,0),\quad e_2^{(78)} = (-6,0,1,0),\quad e_3^{(78)} = (-11,0,0,1). 
\end{eqnarray*}

Define polytopes $\Delta_1,\, \Delta_2,\, \Delta_3,\, \Delta_4$ by 
\begin{eqnarray*}
\Delta_1 & := & \Delta^{(78)} = {\rm Conv}\left\{ \begin{array}{l}(-1,-1,1),\, (-1,-1,-1),\, (4,-1,-1),\\ (3,0,-1),\, (0,2,-1),\, (-1,2,-1)\end{array}\right\}, \\
\Delta_2 & := & {\rm Conv}\left\{ \begin{array}{l}(-1,-1,1),\, (-1,-1,-1),\, (4,-1,-1),\\ (3,0,-1),\, (0,2,-1),\, (-1,1,-1)\end{array}\right\}, \\
\Delta_3 & := & {\rm Conv}\left\{ \begin{array}{l}(-1,-1,1),\, (-1,-1,-1),\, (3,-1,-1),\\ (3,0,-1),\, (0,2,-1),\, (-1,2,-1)\end{array}\right\}, \\
\Delta_4 & := & {\rm Conv}\left\{ \begin{array}{l}(-1,-1,1),\, (-1,-1,-1),\, (3,-1,-1),\\ (3,0,-1),\, (0,2,-1),\, (-1,1,-1)\end{array}\right\}. 
\end{eqnarray*}
Since the polar dual polytopes $\Delta_1^*,\, \Delta_2^*$, and $\Delta_3^*$ of $\Delta_1,\, \Delta_2$, and $\Delta_3$ are the convex hulls of vertices $(1,0,0)$,\, $(0,1,0)$,\, $(0,0,1)$,\, $(-2,-2,-5)$,\, $(-4,-6,-11)$, and $(0,-2,-3)$, resp., $(1,0,0)$,\, $(0,1,0)$,\, $(0,0,1)$,\, $(-2,-2,-5)$,\, $(-4,-6,-11)$, and $(1,-1,-1)$, resp., $(1,0,0)$,\, $(0,1,0)$,\, $(0,0,1)$,\, $(-1,0,-2)$,\, $(-4,-6,-11)$, and $(0,-2,-3)$, and the linear map of $\mathbb{R}^3$ determined by a matrix $\left(\begin{smallmatrix} 1 & 1 & 2\\ 1 & 2 & 3\\ 2 & 3 & 6 \end{smallmatrix}\right)$ gives isomorphisms from $\Delta_4$ to $\Delta_1^*$, from $\Delta_2$ to $\Delta_2^*$, and from $\Delta_3$ to $\Delta_3^*$, the relations $\Delta_1^*\simeq\Delta_4$, $\Delta_2^*\simeq\Delta_2$, and  $\Delta_3^*\simeq\Delta_3$ hold. 

The Newton polytope of $F=F'$, which is the convex hull of vertices  $(-1,-1,1)$,\, $(-1,-1,-1)$,\, $(3,0,-1)$, and $(0,2,-1)$, is not reflexive. 
It is observed that one might take a polytope with vertices $(-1,a,-1)$ with $a$ being $1$ or $2$, and $(b,-1,-1)$ with $b$ being $3$ or $4$ instead of faces spanned by vertices $(-1,-1,1)$,\, $(-1,-1,-1)$,\, $(0,2,-1)$ , and $(-1,-1,1)$,\, $(-1,-1,-1)$,\, $(3,0,-1)$. 
Therefore, there are there polytope-dual pairs, that is, if $(a,b)=(2,4)$ for one side and $(a,b)=(1,3)$ for the other, then, ($\Delta_1, \, \Delta_4)$, if $(a,b)=(1,4)$ for both sides, then, $(\Delta_2,\, \Delta_2)$, if $(a,b)=(2,3)$ for both sides, then, $(\Delta_3, \Delta_3)$. 

\subsection{No. 20} 
We claim that there exist two polytope-dual pairs. 
Take a basis $\{ e_1^{(n)},\, e_2^{(n)},\, e_3^{(n)}\}$ of $M_{n}$ for $n=39, 60$ by 
\begin{eqnarray*}
e_1^{(39)} =  (-3,1,0,0),& & e_2^{(39)} = (-5,0,1,0),\quad e_3^{(39)} = (-9,0,0,1), \\
e_1^{(60)} =  (-4,1,0,0),& & e_2^{(60)} = (-6,0,1,0),\quad e_3^{(60)} = (-7,0,0,1).   
\end{eqnarray*}

Define polytopes $\Delta_1,\, \Delta_1',\, \Delta_2$ and $\Delta_2'$ by 
\begin{eqnarray*}
\Delta_1 & := & {\rm Conv}\left\{\begin{array}{l} (0,-1,1),\, (-1,-1,-1),\\ (3,-1,-1),\, (2,0,-1),\, (-1,2,-1)\end{array}\right\}, \\
\Delta_1' & := & {\rm Conv}\left\{\begin{array}{l} (-1,-1,1),\, (-1,-1,-1),\\ (5,-1,-1),\, (0,2,-1),\, (-1,1,-1)\end{array}\right\}, \\
\Delta_2 & := & {\rm Conv}\left\{\begin{array}{l} (0,-1,1),\, (-1,-1,-1),\\ (2,-1,-1),\, (2,0,-1),\, (-1,2,-1)\end{array}\right\}, \\
\Delta_2' & := & \Delta^{(39)}={\rm Conv}\left\{\begin{array}{l} (-1,-1,1),\, (-1,-1,-1),\\ (5,-1,-1),\, (0,2,-1),\, (-1,2,-1)\end{array}\right\}. 
\end{eqnarray*}
Since the polar dual polytopes $\Delta_1^*$ and $\Delta_2^*$ of $\Delta_1$ and $\Delta_2$ are the convex hulls of vertices $(0,1,0)$,\, $(0,0,1)$,\, $(-2,-2,-3)$,\, $(-4,-6,-7)$, and $(2,0,-1)$, respectively, $(0,1,0)$,\, $(0,0,1)$,\, $(-1,0,-1)$,\, $(-4,-6,-7)$, and $(2,0,-1)$, and the linear map of $\mathbb{R}^3$ determined by a matrix $\left(\begin{smallmatrix} 1 & 1 & 1\\ 1 & 2 & 2\\ 2 & 3 & 4 \end{smallmatrix}\right)$ gives isomorphisms from $\Delta_1'$ to $\Delta_1^*$, and from $\Delta_2$ to $\Delta_2^*$, the relations $\Delta_1^*\simeq\Delta_1'$ and $\Delta_2^*\simeq\Delta_2'$ hold. 

The Newton polytope of $F$, which is the convex hull of vertices $(0,-1,1)$,\, $(-1,-1,-1)$,\, $(2,0,-1)$, and $(-1,2,-1)$, is not reflexive. 
It is observed that one might take a polytope with a vertex $(a,-1,-1)$ with $a$ being $2$ or $3$ instead of a face spanned by vertices $(0,-1,1)$,\, $(-1,-1,-1)$, and $(2,0,-1)$. 
Besides, the Newton polytope of $F'$, which is the convex hull of vertices $(-1,-1,1)$,\, $(-1,-1,-1)$,\, $(2,-1,0)$, and $(0,2,-1)$, is not reflexive. 
It is observed that one might take a polytope with vertices $(-1,b,-1)$ with $b$ being $1$ or $2$, and $(b',-1,-1)$ with $b'$ being $3$ or $5$ instead of faces spanned by vertices $(-1,-1,1)$,\, $(-1,-1,-1)$,\, $(0,2,-1)$ and $(2,-1,0)$,\, $(-1,-1,-1)$,\, $(0,2,-1)$. 
Therefore, there are two polytope-dual  pairs, that is, if $(a,b,b')=(3,1,5)$, then, $(\Delta_1,\, \Delta_1')$, and if $(a,b,b')=(2,2,5)$, then, $(\Delta_2,\, \Delta_2')$. 

\subsection{No. 21} 
We claim that there exist two polytope-dual pairs. 
Take a basis $\{ e_1^{(22)},\, e_2^{(22)},\, e_3^{(22)}\}$ of a lattice $M_{22}$ with  by 
\begin{eqnarray*}
e_1^{(22)} = (-3,1,0,0),\quad e_2^{(22)} = (-5,0,1,0),\quad e_3^{(22)} = (-6,0,0,1). 
\end{eqnarray*}

Define polytopes $\Delta_1,\, \Delta_2$ and $\Delta_2'$ by 
\begin{eqnarray*}
\Delta_1 & := & {\rm Conv}\{ (0,-1,1),\, (-1,-1,-1),\, (4,-1,-1),\, (-1,2,-1)\}, \\
\Delta_2 & := & {\rm Conv}\left\{ \begin{array}{l}(0,-1,1),\, (-1,-1,-1),\\ (2,-1,-1),\, (2,-1,0),\, (-1,2,-1)\end{array}\right\}, \\
\Delta_2' & := & \Delta^{(22)}={\rm Conv}\left\{\begin{array}{l} (-1,-1,1),\, (-1,-1,-1),\\ (4,-1,-1),\, (0,-1,1),\, (-1,2,-1)\end{array}\right\}. 
\end{eqnarray*}
Since the polar dual polytopes $\Delta_1^*$ and $\Delta_2^*$ of $\Delta_1$ and $\Delta_2$ are the convex hulls of vertices $(0,1,0)$,\, $(0,0,1)$,\, $(-3,-5,-6)$, and $(2,0,-1)$, respectively, $(0,1,0)$,\, $(0,0,1)$,\, $(-1,-1,0)$,\, $(-3,-5,-6)$, and $(2,0,-1)$, and the linear map of $\mathbb{R}^3$ determined by a matrix $\left(\begin{smallmatrix} 1 & 1 & 1\\ 1 & 2 & 2\\ 1 & 2 & 3 \end{smallmatrix}\right)$ gives isomorphisms from $\Delta_1$ to $\Delta_1^*$, and from $\Delta_2'$ to $\Delta_2^*$, the relations $\Delta_1^*\simeq\Delta_1$ and $\Delta_2^*\simeq\Delta_2'$ hold. 

The Newton polytope of $F=F'$, which is the convex hull of vertices $(0,-1,1)$,\, $(-1,-1,-1)$,\, $(2,-1,0)$, and $(-1,2,-1)$, is not reflexive. 
It is observed that one might take a polytope with a vertex $(a,-1,-1)$ with $a$ being $2$ or $4$ instead of a face spanned by vertices $(2,-1,0)$,\, $(-1,-1,-1)$, and $(-1,2,-1)$. 
Therefore, there are two polytope-dual pairs, that is, if $a=4$, then, $(\Delta_1, \, \Delta_1)$, and if $a=2$ for one side, and $a=4$ for another, then, $(\Delta_2,\,  \Delta_2')$. 

\subsection{No. 22--No. 24} 
We claim that there exist a unique polytope-dual pair for Nos. 22 and 23, and that none for No. 24. 
Take a basis $\{ e_1^{(n)},\, e_2^{(n)},\, e_3^{(n)}\}$ of lattices $M_{n}$ for $n=9$ and $71$ by 
\begin{eqnarray*}
e_1^{(9)} = (-4,1,0,0), &  & e_2^{(9)} = (-5,0,1,0),\quad e_3^{(9)} = (-10,0,0,1),\, \\
e_1^{(71)} = (-3,1,0,0),& & e_2^{(71)} = (-4,0,1,0),\quad e_3^{(71)} = (-7,0,0,1).   
\end{eqnarray*}

Define polytopes $\Delta_1$ and $\Delta_2$ by 
\begin{eqnarray*}
\Delta_1 & := & \Delta^{(9)}={\rm Conv}\{ (-1,-1,1),\, (-1,-1,-1),\, (4,-1,-1),\, (-1,3,-1)\}, \\
\Delta_2 & := & {\rm Conv}\{ (-1,-1,1),\, (-1,-1,-1),\, (4,-1,-1),\, (-1,1,0)\}. 
\end{eqnarray*}
Since the polar dual polytope $\Delta_1^*$ of $\Delta_1$ is the convex hull of vertices $(1,0,0)$,\, $(0,1,0)$,\, $(0,0,1)$, and $(-4,-5,-10)$, and the linear maps of $\mathbb{R}^3$ determined by a matrix $\left(\begin{smallmatrix} 1 & 1 & 2\\ 1 & 1 & 3\\ 2 & 3 & 5 \end{smallmatrix}\right)$ gives an isomorphism from $\Delta_2$ to $\Delta_1^*$, the relation $\Delta_1^*\simeq\Delta_2$ holds. 

\noindent
{\bf No. 22.}
The Newton polytope of $F=F'$, which is the convex hull of vertices $(-1,-1,1)$,\, $(-1,-1,-1)$,\, $(4,-1,-1)$, and $(-1,1,0)$, coincides with $\Delta_2$. 
Therefore, a pair $(\Delta_1,\, \Delta_{F'}=\Delta_2)$ is polytope-dual. 

\noindent
{\bf No. 23.}
The Newton polytope of $F$, which is the convex hull of vertices $(-1,-1,1)$,\, $(-1,-1,-1)$,\, $(4,-1,-1)$, and $(-1,1,0)$, coincides with $\Delta_2$. 
Besides, the Newton polytope of $F'$, which is the convex hull of vertices $(-1,-1,1)$,\, $(-1,0,-1)$,\, $(4,-1,-1)$, and $(-1,1,0)$, is not reflexive. 
It is observed that one might take a polytope with vertices $(-1,-1,-1)$ and $(-1,3,-1)$ instead of faces spanned by vertices $(-1,-1,1)$,\, $(-1,0,-1)$,\, $(4,-1,-1)$, and $(-1,1,0)$,\, $(-1,0,-1)$,\, $(4,-1,-1)$. 
Therefore, a pair $(\Delta_F=\Delta_2,\, \Delta_1)$ is polytope-dual.  

\noindent
{\bf No. 24.}
The Newton polytope of $F=F'$, which is the convex hull of vertices  $(-1,-1,1)$,\, $(-1,0,-1)$,\, $(4,-1,-1)$, and $(-1,1,0)$, is not reflexive. 
It is observed that one might take a polytope with vertices $(-1,2,-1),\, (0,2,-1)$ and $(-1,-1,-1)$ instead of faces spanned by vertices $(-1,1,0)$,\, $(-1,0,-1)$,\, $(4,-1,-1)$, and $(-1,-1,1)$,\, $(-1,0,-1)$,\, $(4,-1,-1)$, namely, $\Delta^{(71)}$. 

In the polytope $\Delta^{(71)}$, the vertex $(-1, -1, -1)$ is adjacent to three other vertices: vertex $(-1,-1,1)$ with an edge $e_1$, vertex $(-1,2,-1)$ with an edge $e_2$, and vertex $(4,-1,-1)$ with an edge $e_3$. 
On the edges $e_1,\, e_2,\, e_3$, respectively, there are one, two, and four inner lattice points. 
Thus, the polar dual polytope ${\Delta^{(71)}}^*$ must contain a triangle as a two-dimensional face that is adjacent to other two-dimensional faces with inner lattice points one, two, and four. 
However, it is easily observed that there is no such a configuration in the polytope ${\Delta^{(71)}}^*$. 
Thus, $\Delta^{(71)}$ is not self-dual, and there does not exist a polytope-dual pair.

\subsection{No. 25} 
We claim that there exist two polytope-dual pairs. 
Take a basis $\{ e_1^{(n)},\, e_2^{(n)},\, e_3^{(n)}\}$ of lattices $M_{n}$ for $n=37, 58$ by 
\begin{eqnarray*}
e_1^{(37)} = (-3,1,0,0), & & e_2^{(37)} = (-4,0,1,0),\quad e_3^{(37)} = (-8,0,0,1), \\
e_1^{(58)} = (-4,1,0,0),& & e_2^{(58)} = (-5,0,1,0),\quad e_3^{(58)} = (-6,0,0,1). 
\end{eqnarray*}

Define polytopes $\Delta_1,\, \Delta_1',\, \Delta_2$ and $\Delta_2'$ by 
\begin{eqnarray*}
\Delta_1 & := & {\rm Conv}\left\{\begin{array}{l} (-1,1,0),\, (-1,-1,0),\\ (-1,-1,-1),\, (3,-1,-1),\, (0,-1,1)\end{array}\right\}, \\
\Delta_1' & := & \Delta^{(37)}={\rm Conv}\left\{\begin{array}{l} (-1,-1,1),\, (-1,-1,-1),\\ (4,-1,-1),\, (3,0,-1),\, (-1,3,-1)\end{array}\right\}, \\
\Delta_2 & := & {\rm Conv}\left\{\begin{array}{l} (-1,1,0),\, (-1,-1,1),\\ (-1,-1,-1),\, (3,-1,-1),\, (0,-1,1)\end{array}\right\}, \\
\Delta_2' & := & {\rm Conv}\left\{\begin{array}{l} (-1,-1,1),\, (-1,-1,-1),\\ (3,-1,-1),\, (3,0,-1),\, (-1,3,-1)\end{array}\right\}. 
\end{eqnarray*}
Since the polar dual polytopes $\Delta_1^*$ and $\Delta_2^*$ of $\Delta_1$ and $\Delta_2$ are the convex hulls of vertices $(1,0,0)$,\, $(0,1,0)$,\, $(0,-1,2)$,\, $(1,0,-1)$, and $(-4, -5, -6)$, respectively, $(1,0,0),$\, $(0,1,0)$,\, $(0,-1,2)$,\, $(0,-1,-2)$, and $(-4, -5, -6)$, and the linear map of $\mathbb{R}^3$ determined by a matrix $\left(\begin{smallmatrix} 1 & 1 & 1\\ 1 & 1 & 2\\ 2 & 3 & 3 \end{smallmatrix}\right)$ gives isomorphisms from $\Delta_1'$ to $\Delta_1^*$, and from $\Delta_2'$ to $\Delta_2^*$, the relations  $\Delta_1^*\simeq\Delta_1'$ and $\Delta_2^*\simeq\Delta_2'$ hold. 

The Newton polytopes of $F$, which is the convex hull of vertices $(0,-1,1)$,\, $(-1,-1,-1)$,\, $(3,-1,-1)$, and $(-1,1,0)$, is not reflexive. 
It is observed that one might take a polytope with vertices $(-1,-1,a)$ with $a$ being $0$ or $1$, and $(-1,a',-1)$ with $a'$ being $1$ or $2$ instead of a face spanned by vertices $(0,-1,1)$,\, $(-1,-1,-1)$, and $(-1,1,0)$. 
Besides, the Newton polytopes of $F'$, which is the convex hull of vertices $(-1,-1,1)$,\, $(-1,-1,-1)$,\, $(3,0,-1)$, and $(-1,1,0)$, is not reflexive. 
It is observed that one might take a polytope with vertices $(b,-1,-1)$ with $b$ being $3$ or $4$, and $(-1,b',-1)$ with $b'$ being $0,\, 2$ or $3$ instead of faces spanned by vertices $(-1,-1,1)$,\, $(-1,-1,-1)$,\, $(3,0,-1)$, and $(-1,1,0)$,\, $(-1,-1,-1)$,\, $(3,0,-1)$. 
Therefore, there are two polytope-dual pairs, that is, if $(a,b,b')=(0, 4, 3)$, then, $(\Delta_1,\, \Delta_1')$, and if $(a,b,b')=(1, 3, 3)$, then, $(\Delta_2,\, \Delta_2')$. 

\subsection{No. 26} 
We claim that there exist four polytope-dual pairs. 
Take a basis $\{ e_1^{(87)},\, e_2^{(87)},\, e_3^{(87)}\}$ of a lattice $M_{87}$ by 
\begin{eqnarray*}
e_1^{(87)} = (-3,1,0,0),\quad e_2^{(37)} = (-4,0,1,0),\quad e_3^{(37)} = (-5,0,0,1). 
\end{eqnarray*}

Define polytopes $\Delta_i,\, \Delta_i'$ for $i=1,2,3,4$ by 
\begin{eqnarray*}
\Delta_1 & := & {\rm Conv}\left\{ \begin{array}{l}(-1,-1,1),\, (-1,1,0),\, (0,-1,1),\\ (-1,-1,-1),\, (2,-1,-1),\, (2,0,-1),\, (-1,0,-1)\end{array}\right\}, \\
\Delta_1' & := & {\rm Conv}\left\{\begin{array}{l} (-1,-1,1),\, (-1,1,0),\, (0,-1,1),\\ (-1,-1,-1),\, (2,-1,-1),\, (2,0,-1),\, (-1,2,-1)\end{array}\right\}, \\
\Delta_2 & := & {\rm Conv}\left\{ \begin{array}{l}(-1,-1,0),\, (-1,1,0),\, (0,-1,1),\\ (-1,-1,-1),\, (3,-1,-1),\, (2,0,-1),\, (-1,0,-1)\end{array}\right\}, \\
\Delta_2' & := & {\rm Conv}\left\{\begin{array}{l} (-1,-1,0),\, (-1,1,0),\, (0,-1,1),\\ (-1,-1,-1),\, (3,-1,-1),\, (2,0,-1),\, (-1,2,-1)\end{array}\right\}, \\
\Delta_3 & := & {\rm Conv}\left\{ \begin{array}{l}(-1,-1,0),\, (-1,1,0),\, (0,-1,1),\\ (-1,-1,-1),\, (2,-1,-1),\, (2,0,-1),\, (-1,2,-1)\end{array}\right\}, \\
\Delta_3' & := & {\rm Conv}\left\{\begin{array}{l} (-1,-1,1),\, (-1,1,0),\, (0,-1,1),\\ (-1,-1,-1),\, (3,-1,-1),\, (2,0,-1),\, (-1,0,-1)\end{array}\right\}, \\
\Delta_4 & := & {\rm Conv}\left\{ \begin{array}{l}(-1,-1,0),\, (-1,1,0),\, (0,-1,1),\\ (-1,-1,-1),\, (2,-1,-1),\, (2,0,-1),\, (-1,0,-1)\end{array}\right\}, \\
\Delta_4' & := & \Delta^{(87)} = {\rm Conv}\left\{\begin{array}{l} (-1,-1,1),\, (-1,1,0),\, (0,-1,1),\\ (-1,-1,-1),\, (3,-1,-1),\, (2,0,-1),\, (-1,2,-1)\end{array}\right\}. 
\end{eqnarray*}
Since the polar dual polytopes $\Delta_i^*,\, i=1,2,3,4$ are respectively given as follows: 
\begin{eqnarray*}
\Delta_1^* & := & {\rm Conv}\left\{ \begin{array}{l}(1,0,0),\, (0,1,0),\, (0,0,1),\, (0,-1,-2),\\ (-1,0,-1),\, (0,-1,1),\, (-3,-4,-5)\end{array}\right\}, \\
\Delta_2^* & := & {\rm Conv}\left\{ \begin{array}{l}(1,0,0),\, (0,1,0),\, (0,0,1),\, (1,0,-1),\\ (-2,-2,-3),\, (0,-1,1),\, (-3,-4,-5)\end{array}\right\}, \\
\Delta_3^* & := & {\rm Conv}\left\{ \begin{array}{l}(1,0,0),\, (0,1,0),\, (0,0,1),\, (1,0,-1),\\ (-1,0,-1),\, (-2,-3,-3),\, (-3,-4,-5)\end{array}\right\}, \\
\Delta_4^* & := & {\rm Conv}\left\{ \begin{array}{l}(1,0,0),\, (0,1,0),\, (0,0,1),\, (1,0,-1),\\ (-1,0,-1),\, (0,-1,1),\, (-3,-4,-5)\end{array}\right\}, 
\end{eqnarray*}
and the linear map of $\mathbb{R}^3$ determined by a matrix $\left(\begin{smallmatrix} 1 & 1 & 1\\ 1 & 1 & 2\\ 1 & 2 & 2 \end{smallmatrix}\right)$ gives isomorphisms from $\Delta_i'$ to $\Delta_i^*,\, i=1,2,3,4$, the relations $\Delta_i^*\simeq\Delta_i',\, i=1,2,3,4$ hold. 

The Newton polytope of $F=F'$, which is the convex hull of vertices  $(0,-1,1)$,\, $(-1,-1,-1)$,\, $(2,0,-1)$, and $(-1,1,0)$, is not reflexive. 
It is observed that one might take a polytope with vertices $(-1,-1,a)$ with $a$ being $0$ or $1$, and $(b,-1,-1)$, with $b$ being $2$ or $3$, and $(-1,c,-1)$ with $c$ being $0$ or $2$  instead of faces spanned by $(-1,1,0)$,\, $(-1,-1,-1)$,\, $(0,-1,1)$, and $(2,0,-1)$,\, $(-1,-1,-1)$,\, $(0,-1,1)$, and $(-1,1,0)$,\, $(-1,-1,-1)$,\, $(2,0,-1)$. 
Therefore, there are four polytope-dual pairs, that is, if $(a,b,c)=(1, 2, 0)$ for one side, and $(a,b,c)=(1, 2, 2)$ for the other, then, $(\Delta_1,\, \Delta_1')$; if $(a,b,c)=(0, 3, 0)$ for one side, and $(a,b,c)=(0, 3, 2)$ for the other, then, $(\Delta_2,\, \Delta_2')$; if $(a,b,c)=(0, 2, 2)$ for one side, and $(a,b,c)=(1, 3, 0)$ for the other, then, $(\Delta_3,\, \Delta_3')$; and if $(a,b,c)=(0, 2, 0)$ for one side, and $(a,b,c)=(1, 3, 2)$ for the other, then, $(\Delta_4,\, \Delta_4')$. 

\subsection{No. 27} 
We claim that there does not exist a polytope-dual pair. 
Take a basis $\{ e_1^{(4)},\, e_2^{(4)},\, e_3^{(4)}\}$ of a lattice $M_{4}$ by 
\begin{eqnarray*}
e_1^{(4)} = (-3,1,0,0),\quad e_2^{(4)} = (-4,0,1,0),\quad e_3^{(4)} = (-4,0,0,1). 
\end{eqnarray*}
By a direct computation, one has ${\Delta^{(4)}}^*\not\simeq\Delta^{(4)}$, which we say {\it not self-dual}. 

It is observed that there exist four invertible projectivisations
\begin{eqnarray*}
F_1 = F_1' & = & X^4+Y^3+Z^3+W^{12}, \\
F_2 =F_2' = X^4+Y^3+Z^3+W^{8}Z, & & F_3 = F_3'= X^4+Y^3+Z^3+W^{8}Y, \\
F_4 = F_4' & = & X^4+Y^3+Z^3+W^{9}X. 
\end{eqnarray*}

The Newton polytope of $F_1=F_1'$, which is the convex hull of vertices $(-1,-1,2)$,\, $(-1,-1,-1)$,\, $(3,-1,-1)$, and $(-1,2,-1)$, coincide with  $\Delta^{(4)}$. 

Since $F_2$ and $F_3$ are symmetric in variables $Y$ and $Z$, we only treat with $F_2$. 
The Newton polytope of $F_2=F_2'$, which is the convex hull of vertices $(-1,-1,2)$,\, $(-1,-1,0)$,\, $(3,-1,-1)$, and $(-1,2,-1)$, is not reflexive by a direct computation. 
It is observed that one might take a polytope with a vertiex $(-1,-1,-1)$ instead of a face spanned by vertices $(-1,-1,0)$,\, $(3,-1,-1)$, and $(-1,2,-1)$, which is $\Delta^{(4)}$. 

The Newton polytope of $F_4=F_4'$, which is the convex hull of vertices $(-1,-1,2)$,\, $(0,-1,-1)$,\, $(3,-1,-1)$, and $(-1,2,-1)$, is not reflexive. 
It is observed that one must take the polytope $\Delta^{(4)}$. 

\subsection{No. 28--No. 29} 
We claim that there exist two polytope-dual pairs for Nos. 28 and 29. 
Take a basis $\{ e_1^{(n)},\, e_2^{(n)},\, e_3^{(n)}\}$ of lattices $M_{n}$ for $n=12, 27$ and $49$ by 
\begin{eqnarray*}
e_1^{(12)} = (-2,1,0,0), & & e_2^{(12)} = (-6,0,1,0),\quad e_3^{(12)} = (-9,0,0,1),\, \\
e_1^{(27)} = (-1,7,-1,-1),& & e_2^{(27)} = (-1,-1,2,-1),\quad e_3^{(27)} = (0,-1,-1,1),\, \\
e_1^{(49)} = (0,7,-1,-1),& & e_2^{(49)} = (-1,-1,2,-1),\quad e_3^{(49)} = (-1,-1,-1,1).   
\end{eqnarray*}

In \cite{KM}, it is proved that the polytopes $\Delta^{(27)}$ and $\Delta^{(49)}$ are isomorphic to the polytope $\Delta_{(27,49)}$ with vertices $(0,0,1)$,\, $(1,0,0)$,\, $(0,1,0)$, and $(-3,-8,-12)$. 

Define polytopes $\Delta_1,\, \Delta_1',\, \Delta_2$ and $\Delta_2'$ by 
\begin{eqnarray*}
\Delta_1 & := & {\rm Conv}\{ (0,0,1),\, (1,0,0),\, (0,1,0),\, (-2,-6,-9)\}, \\
\Delta_1' & := & \Delta^{(12)}={\rm Conv}\{ (-1,-1,1),\, (-1,-1,-1),\, (8,-1,-1),\, (-1,2,-1)\}, \\
\Delta_2 & := & \Delta^{(49)}={\rm Conv}\{ (0,0,1),\, (1,0,0),\, (0,1,0),\, (-3,-8,-12)\}, \\
\Delta_2' & := & {\rm Conv}\{ (-1,-1,1),\, (-1,-1,-1),\, (7,-1,-1),\, (-1,2,-1)\}. 
\end{eqnarray*}
It is straight-forward to see that the relations $\Delta_1^*\simeq\Delta_1' $ and $\Delta_2^*\simeq\Delta_2' $ hold. 

\noindent
{\bf No. 28.}
The Newton polytope of $F$, which is the convex hull of vertices $(1,0,0)$,\, $(0,1,0)$,\, $(0,0,1)$, and $(-2,-6,-9)$, coincides with $\Delta_1'$. 
Besides, the Newton polytope of $F'$, which is the convex hull of vertices $(-1,-1,1)$,\, $(-1,-1,0)$,\, $(7,-1,-1)$, and $(-1,2,-1)$, is not reflexive. 
It is observed that one might take a polytope with a vertex $(-1,-1,-1)$ instead of a face spanned by vertices $(-1,-1,0)$,\, $(7,-1,-1)$, and $(-1,2,-1)$. 
Therefore, there are two polytope-dual pairs, that is, $(\Delta_{F} = \Delta_1,\, \Delta_1')$, and $(\Delta_2,\, \Delta_2')$. 

\noindent
{\bf No. 29.}
The Newton polytope of $F$, which is the convex hull of vertices $(1,0,0)$,\, $(0,1,0)$,\, $(0,0,1)$, and $(-2,-6,-9)$, coincides with $\Delta_1'$. 
Besides, the Newton polytope of $F'$, which is the convex hull of vertices $(-1,-1,1)$,\, $(0,-1,-1)$,\, $(7,-1,-1)$, and $(-1,2,-1)$, is not reflexive. 
It is observed that one might take a polytope with vertices $(-1,-1,-1)$ and $(-1,0,-1)$ instead of a face spanned by vertices $(-1,-1,1)$,\, $(0,-1,-1)$, and $(-1,2,-1)$. 
Therefore, there are two polytope-dual pairs, that is, $(\Delta_{F} = \Delta_1,\, \Delta_1')$, and $(\Delta_2,\, \Delta_2')$. 

\subsection{No. 30} 
We claim that there exist two polytope-dual pairs. 
Take a basis $\{ e_1^{(n)},\, e_2^{(n)},\, e_3^{(n)}\}$ of lattices $M_{n}$ for $n=40, 81$  by 
\begin{eqnarray*}
e_1^{(40)} = (-2,1,0,0),& & e_2^{(40)} = (-4,0,1,0),\quad e_3^{(40)} = (-7,0,0,1), \\
e_1^{(81)} = (0,7,-1,-1),& & e_2^{(81)} = (0,-1,2,-1),\quad e_3^{(81)} = (-1,-1,-1,1). 
\end{eqnarray*}

Define polytopes $\Delta_1,\, \Delta_1',\, \Delta_2,\, \Delta_2'$ by 
\begin{eqnarray*}
\Delta_1 & := & {\rm Conv}\left\{\begin{array}{l} (0,1,0),\, (0,0,1),\\ (1,1,1),\, (0,-2,-3),\, (-2,-6,-9)\end{array}\right\}, \\
\Delta_1' & := & \Delta^{(40)}={\rm Conv}\left\{\begin{array}{l} (-1,-1,1),\, (-1,-1,-1),\\ (6,-1,-1),\, (0,2,-1),\, (-1,2,-1)\end{array}\right\}, \\
\Delta_2 & := & {\rm Conv}\left\{\begin{array}{l} (0,1,0),\, (0,0,1),\\ (1,1,1),\, (1,0,0),\, (-2,-6,-9)\end{array}\right\}, \\
\Delta_2' & := & {\rm Conv}\left\{\begin{array}{l} (-1,-1,1),\, (-1,-1,-1),\\ (6,-1,-1),\, (0,2,-1),\, (-1,1,-1)\end{array}\right\}. 
\end{eqnarray*}
Since the polar dual polytopes $\Delta_1^*$ and $\Delta_2^*$ of $\Delta_1$ and $\Delta_2$ are respectively the convex hulls of vertices $(-1,-1,1)$,\, $(1,-1,-1)$,\, $(-2,2,-1)$,\, $(-1,2,-1)$, and $(8,-1,-1)$, respectively, $(-1,-1,1)$,\, $(1,-1,-1)$,\, $(-1,1,-1)$,\, $(-1,2,-1)$, and $(8,-1,-1)$, and the linear map of $\mathbb{R}^3$ determined by a matrix $\left(\begin{smallmatrix} 1 & 0 & 0\\ -1 & 1 & 0\\ -1 & 0 & 1 \end{smallmatrix}\right)$ gives isomorphisms from $\Delta_i'$ to $\Delta_i^*$ for $i=1,2$, the relations $\Delta_1^*\simeq\Delta_1'$ and $\Delta_2^*\simeq\Delta_2'$ hold. 

The Newton polytope of $F$, which is the convex hull of vertices $(1,1,1)$,\, $(0,1,0)$,\, $(0,0,1)$, and $(-2,-6,-9)$, is not reflexive. 
It is observed that one might take a polytope with vertices $(1,0,0)$ and $(0,-2,-3)$ instead of a face spanned by vertices $(1,1,1)$,\, $(0,0,1)$, and $(-2,-6,-9)$. 
Besides, the Newton polytope of $F'$, which is the convex hull of vertices $(-1,-1,1)$,\, $(-1,0,-1)$,\, $(5,-1,-1)$, and $(0,2,-1)$, is not reflexive. 
It is observed that one might take a polytope with vertices $(-1,-1,-1)$ and $(-1,b,-1)$ with $b$ being $1$ or $2$ instead of faces spanned by vertices $ (-1,-1,1)$,\, $(-1,0,-1)$,\, $(0,2,-1)$, and $(-1,-1,1)$,\, $(-1,0,-1)$,\, $(5,-1,-1)$. 
Therefore, there are two polytope-dual pairs, that is, $(\Delta_1, \, \Delta_1')$, and $(\Delta_2,\, \Delta_2')$. 

\subsection{No. 31} 
We claim that there exist four polytope-dual pairs. 
Take a basis $\{ e_1^{(n)},\, e_2^{(n)},\, e_3^{(n)}\}$ of lattices $M_{n}$ for $n=11, 24$ by 
\begin{eqnarray*}
e_1^{(11)} = (2,7,-1,-1),& & e_2^{(11)} = (-1,-1,2,-1),\quad e_3^{(11)} = (-1,-1,-1,1), \\
e_1^{(24)} = (-2,1,0,0),& & e_2^{(24)} = (-4,0,1,0),\quad e_3^{(24)} = (-5,0,0,1). 
\end{eqnarray*}

Define polytopes $\Delta_i,\, \Delta_i',\, i=1,2,3,4$ by 
\begin{eqnarray*}
\Delta_1 & := & \Delta^{(11)}={\rm Conv}\{ (0,1,0),\, (0,0,1),\, (2,2,3),\, (-3,-8,-12)\}, \\
\Delta_1' & := & {\rm Conv}\{ (0,-1,1),\, (-1,-1,-1),\, (4,-1,-1),\, (-1,2,-1)\}, \\
\Delta_2 & := & {\rm Conv}\{ (0,1,0),\, (0,0,1),\, (2,2,3),\, (-2,-6,-9)\}, \\
\Delta_2' & := & {\rm Conv}\{ (0,-1,1),\, (-1,-1,-1),\, (5,-1,-1),\, (-1,2,-1)\}, \\
\Delta_3 & := & {\rm Conv}\{ (0,1,0),\, (0,0,1),\, (1,1,2),\, (0,-2,-3),\, (-2,-6,-9)\}, \\
\Delta_3' & := & \Delta^{(24)}={\rm Conv}\left\{\begin{array}{l} (0,-1,1),\, (-1,-1,-1),\\ (5,-1,-1),\, (-1,-1,1),\, (-1,2,-1)\end{array}\right\}, \\
\Delta_4 & := & {\rm Conv}\left\{ \begin{array}{l}(0,1,0),\, (0,0,1),\, (1,1,2),\\ (0,-2,-3),\, (-2,-6,-9),\, (-1,-2,-4)\end{array}\right\}, \\
\Delta_4' & := & {\rm Conv}\left\{\begin{array}{l} (0,-1,1),\, (-1,-1,-1),\, (4,-1,-1),\\ (3,0,-1),\, (-1,2,-1),\, (-1,-1,1)\end{array}\right\}. 
\end{eqnarray*}
Since the polar dual polytopes $\Delta_i^*$ of $\Delta_i,\, i=1,2,3,4$ are respectively the convex hulls of vertices $(-1,-1,1)$,\, $(2,-1,-1)$,\, $(-1,2,-1)$, and $(7,-1,-1)$, resp., $(-1,-1,1)$,\, $(2,-1,-1)$,\, $(-1,2,-1)$, and $(8,-1,-1)$, resp., $(-1,-1,1)$,\, $(2,-1,-1)$,\, $(-1,2,-1)$,\, $(-2,-1,1)$, and $(8,-1,-1)$, resp., $(-1,-1,1)$,\, $(2,-1,-1)$,\, $(-1,2,-1)$,\, $(-2,-1,1)$,\, $(5,0,-1)$, and $(7,-1,-1)$, and the linear map of $\mathbb{R}^3$ determined by a matrix $\left(\begin{smallmatrix} 1 & 1 & 1\\ 1 & 2 & 2\\ 1 & 2 & 3 \end{smallmatrix}\right)$ gives isomorphisms from $\Delta_i'$ to $\Delta_i^*$ for $i=1,2,3,4$, the relations $\Delta_i^*\simeq\Delta_i',\, i=1,2,3,4$ hold. 

The Newton polytope of $F$, which is the convex hull of vertices $(0,1,0)$,\, $(0,0,1)$,\, $(1,1,2)$,\, $(-2,-6,-9)$, is not reflexive. 
It is observed that one might take a polytope with a vertex $(0,-2,-3)$, $(2,2,3)$, $(-1,-2,-4)$, $(-2,-6,-2)$, or $(-3,-8,-12)$ instead of a face spanned by vertices $(1,1,2)$,\, $(0,1,0)$, and $(-2,-6,-9)$. 
Besides, the Newton polytope of $F'$, which is the convex hull of vertices $(0,-1,1)$,\, $(-1,-1,-1)$,\, $(4,-1,-1)$, and $(-1,2,-1)$, coincides with $\Delta_1'$. 
Therefore, there are four polytope-dual pairs, that is, $(\Delta_1,\, \Delta_{F'}=\Delta_1')$, $(\Delta_2,\, \Delta_2')$, $(\Delta_3,\, \Delta_3')$, and $(\Delta_4,\, \Delta_4')$. 

\subsection{No. 32--No. 34} 
We claim that there exist two polytope-dual pairs for Nos. 32 to 34. 
Take a basis $\{ e_1^{(n)},\, e_2^{(n)},\, e_3^{(n)}\}$ of lattices $M_{n}$ with  for $n=6,26,34,$ and $76$ by 
\begin{eqnarray*}
e_1^{(6)} = (-2,1,0,0),& & e_2^{(6)} = (-2,0,1,0),\quad e_3^{(6)} = (-5,0,0,1),\, \\
e_1^{(26)} = (-1,4,-1,-1),& & e_2^{(26)} = (-1,-1,3,-1),\quad e_3^{(26)} = (0,-1,-1,1),\, \\ 
e_1^{(34)} = (-1,4,-1,-1),& & e_2^{(34)} = (0,-1,3,-1),\quad e_3^{(34)} = (-1,-1,-1,1),\, \\
e_1^{(76)} = (1,1,1,-1),& & e_2^{(76)} = (4,1,0,-1),\quad e_3^{(76)} = (9,-1,0,-1).   
\end{eqnarray*}

In \cite{KM}, it is proved that the polytopes $\Delta^{(26)}$ and $\Delta^{(34)}$ are isomorphic to the polytope with vertices $(1,0,0)$,\, $(0,1,0)$,\, $(0,0,1)$, and $(-4,-5,-10)$. 

Define polytopes $\Delta_1,\, \tilde{\Delta}_1,\, \Delta_1'$, and $\Delta_2,\, \tilde{\Delta}_2,\, \Delta_2'$ by 
\begin{eqnarray*}
\Delta_1 & := & {\rm Conv}\{  (1,0,0),\, (0,1,0),\, (0,0,1),\, (-2,-2,-5)\}, \\
\tilde{\Delta}_1 & := & {\rm Conv}\{  (0,0,-1),\, (2,0,-1),\, (2,-1,0),\, (-3,1,3)\}, \\
\Delta_1' & := & \Delta^{(6)}={\rm Conv}\{ (-1,-1,1),\, (-1,-1,-1),\, (4,-1,-1),\, (-1,4,-1)\}, \\
\Delta_2 & := & {\rm Conv}\left\{\begin{array}{l}  (1,0,0),\, (0,1,0),\, (0,0,1),\\ (-2,-2,-5),\, (-a, -(a+1), -2(a+1))\end{array}\right\}, \\
\tilde{\Delta}_2 & := & {\rm Conv}\left\{\begin{array}{l}  (0,0,-1),\, (2,0,-1),\, (2,-1,0),\\ (-3,1,3),\, (a-2, 1, -a+2)\end{array}\right\}, \\
\Delta_2' & := & {\rm Conv}\left\{\begin{array}{l} (-1,-1,1),\, (-1,-1,-1),\, (4,-1,-1),\\ (-1,3,-1),\, (a, 3-a,-1)\end{array}\right\}, \\
 & & \textnormal{here}, \quad a = 3,2,1,0. 
\end{eqnarray*}
Since the polar dual polytopes $\Delta_1^*$ and $\tilde{\Delta}_1^*$ of $\Delta_1$ and $\tilde{\Delta}_1$ are the convex hulls of vertices $(-1,-1,1)$,\, $(-1,-1,-1)$,\, $(4,-1,-1)$, and $(-1,4,-1)$, respectively, $(-1,-1,-1)$,\, $(0,1,1)$,\, $(0,-4,1)$, and $(5,11,1)$, and the linear map of $\mathbb{R}^3$ determined by a matrix $\left(\begin{smallmatrix} -1 & -3 & 0\\ -1 & -2 & 0\\ -3 & -6 & -1 \end{smallmatrix}\right)$ gives an isomorphism from $\Delta_1'$ to $\tilde{\Delta}_1^*$, the relations $\Delta_1^*=\Delta_1'\simeq\tilde{\Delta}_1^*$ hold. 
Moreover, Since the polar dual polytopes $\Delta_2^*$ and $\tilde{\Delta}_2^*$ of $\Delta_2$ and $\tilde{\Delta}_2$ are convex hulls of the vertices $(-1,-1,1)$,\, $(-1,-1,-1)$,\, $(4,-1,-1)$,\,  $(a, 3-a, -1)$,and $(-1,3,-1)$, respectively, $(-1,-1,-1)$,\, $(0,1,1)$,\, $(5,11,1)$,\, $(1,-1,1)$, and $(0, a-3, 1)$, and the linear map of $\mathbb{R}^3$ determined by a matrix $\left(\begin{smallmatrix} -1 & -2 & 0\\ -1 & -3 & 0\\ -3 & -6 & -1 \end{smallmatrix}\right)$ gives an isomorphism from $\Delta_2'$ to $\tilde{\Delta}_2^*$, the relations $\Delta_2^* =\Delta_2' \simeq\tilde{\Delta}_2^*$ hold. 

\noindent
{\bf No. 32.}
The Newton polytope of $F$, which is the convex hull of vertices $(1,0,0)$,\, $(0,1,0)$,\, $(0,0,1)$, and $(-2,-2,-5)$, coincides with $\Delta_1$. 
Besides, the Newton polytope of $F'$, which is the convex hull of vertices $(-1,-1,0)$,\, $(-1,3,-1)$,\, $(-1,-1,1)$, and $(4,-1,-1)$, is not reflexive. 
It is observed that one might take a polytope with a vertex $(-1,-1,-1)$ instead of a face spanned by vertices $(-1,-1,0)$,\, $(-1,-1,1)$, and $(-1,3,-1)$. 
Therefore, there are two polytope-dual pairs, that is, $(\Delta_{F} = \Delta_1,\, \Delta_1')$, and $(\Delta_2,\, \Delta_2')$. 

\noindent
{\bf No. 33.}
The Newton polytope of $F$, which is the convex hull of vertices $(1,0,0)$,\, $(0,1,0)$,\, $(0,0,1)$, and $(-2,-2,-5)$, coincides with $\Delta_1$. 
Besides, the Newton polytope of $F'$, which is the convex hull of vertices $(-1,0,-1)$,\, $(-1,3,-1)$,\, $(-1,-1,1)$, and $(4,-1,-1)$, is not reflexive. 
It is observed that one might take a polytope with a vertex $(-1,-1,-1)$ instead of a face spanned by vertices $(-1,-1,1)$,\, $(-1,0,-1)$, and $(4,-1,-1)$. 
Therefore, there are two polytope-dual pairs, that is, $(\Delta_{F} = \Delta_1,\, \Delta_1')$, and $(\Delta_2,\, \Delta_2')$. 

\noindent
{\bf No. 34.}
The Newton polytope of $F$, which is the convex hull of vertices $(0,0,-1)$,\, $(2,0,-1)$,\, $(2,-1,0)$, and $(-3,1,3)$, coincides with $\tilde{\Delta}_1$. 
Besides, the Newton polytope of $F'$, which is the convex hull of vertices $(-1,0,-1)$,\, $(3,-1,-1)$,\, $(0,3,-1)$, and $(-1,-1,1)$, is not reflexive. 
It is observed that one might take a polytope with vertices $(-1,-1,-1)$, and $(-1,b',-1)$ with $b'$ being $3$ or $4$ instead of faces spanned by vertices $(-1,-1,1)$,\, $(-1,0,-1)$,\, $(3,-1,-1)$, and $(-1,-1,1)$,\, $(-1,0,-1)$, \, $(0,3,-1)$. 
Therefore, there are two polytope-dual pairs, that is, $(\Delta_F=\tilde{\Delta}_1,\, \Delta_1')$, and $(\tilde{\Delta}_2,\, \Delta_2')$. 

\subsection{No. 35--No. 37} 
We claim that there exists a unique polytope-dual pair for Nos. 35 to 37. 
Take a basis $\{ e_1^{(n)},\, e_2^{(n)},\, e_3^{(n)}\}$ lattices $M_{n}$ for $n=10,46,65,$ and $80$ by 
\begin{eqnarray*}
e_1^{(10)} = (-1,1,0,0),& & e_2^{(10)} = (-4,0,1,0),\quad e_3^{(10)} = (-6,0,0,1),\, \\
e_1^{(46)} = (1,1,1,-1),& & e_2^{(46)} = (2,2,-1,0),\quad e_3^{(46)} = (11,0,-1,-1),\, \\ 
e_1^{(65)} = (1,0,1,-1),& & e_2^{(65)} = (2,1,-1,0),\quad e_3^{(65)} = (10,-1,-1,-1),\, \\
e_1^{(80)} = (1,1,1,-1),& & e_2^{(80)} = (2,1,-1,0),\quad e_3^{(80)} = (10,-1,-1,-1).   
\end{eqnarray*}

In \cite{KM}, it is proved that the polytopes $\Delta^{(46)},\, \Delta^{(65)}$ and $\Delta^{(80)}$ are isomorphic to the polytope $\Delta_{(46,65,80)}$ with vertices $(-1,0,0)$,\, $(1,-1,0)$,\, $(0,0,1)$, and $(2,4,-1)$. 
Let a polytope $\Delta_1$ be $\Delta_{(46,65,80)}$ and $\Delta_1'$ be the polytope $\Delta^{(10)}$ which is the convex hull of vertices $(-1,-1,1)$,\, $(-1,-1,-1)$,\, $(11,-1,-1)$, and $(-1,2,-1)$.
Since the polar dual polytope $\Delta_1^*$ of $\Delta_1$ is the convex hull of vertices $(1,-1,-1)$,\, $(1,2,11)$,\, $(1,2,-1)$, and $(-1,0,-1)$, and the linear map of $\mathbb{R}^3$ determined by a matrix $\left(\begin{smallmatrix} 0 & 0 & -1\\ 0 & -1 & -4\\ -1 & -1 & -6 \end{smallmatrix}\right)$ gives an isomorphism from $\Delta_1'$ to $\Delta_1^*$, the relation $\Delta_1^*\simeq\Delta_1'$ holds. 

\noindent
{\bf No. 35.}
The Newton polytope of $F$, which is the convex hull of vertices $(-1,0,0)$,\, $(1,-1,0)$,\, $(0,0,1)$, and $(2,4,-1)$, coincides with $\Delta_1$. 
Besides, the Newton polytope of $F'$, which is the convex hull of vertices $(-1,-1,1)$,\, $(0,-1,-1)$,\, $(5,-1,0)$, and $(-1,2,-1)$, is not reflexive. 
It is observed that the pair $(\Delta_{F} = \Delta_1,\, \Delta_1')$ is polytope-dual. 

\noindent
{\bf No. 36.}
The Newton polytope of $F$, which is the convex hull of vertices $(-1,0,0)$,\, $(1,-1,0)$,\, $(0,0,1)$, and $(2,4,-1)$, coincides with $\Delta_1$. 
Besides, the Newton polytope of $F'$, which is the convex hull of vertices $(-1,-1,1)$,\, $(0,-1,-1)$,\, $(7,0,-1)$, and $(-1,2,-1)$, is not reflexive. 
It is observed that the pair $(\Delta_{F} = \Delta_1,\, \Delta_1')$ is polytope-dual. 

\noindent
{\bf No. 37.}
The Newton polytope of $F$, which is the convex hull of vertices $(-1,0,0)$,\, $(1,-1,0)$,\, $(0,0,1)$, and $(2,4,-1)$, coincides with $\Delta_1$. 
Besides, the Newton polytope of $F'$, which is the convex hull of vertices $(-1,-1,1)$,\, $(-1,-1,-1)$,\, $(10,-1,-1)$, and $(-1,2,-1)$, is not reflexive. 
It is observed that the pair $(\Delta_{F} = \Delta_1,\, \Delta_1')$ is polytope-dual. 

\subsection{No. 38--No. 40} 
We claim that there exist two polytope-dual pairs for Nos. 38 and 40, and that none for No. 39. 
Take a basis $\{ e_1^{(n)},\, e_2^{(n)},\, e_3^{(n)}\}$ of lattices $M_{n}$ for $n=42,68,83,$ and $92$ by 
\begin{eqnarray*}
e_1^{(42)} = (-1,1,0,0),& & e_2^{(42)} = (-3,0,1,0),\quad e_3^{(42)} = (-5,0,0,1),\, \\
e_1^{(68)} = (1,0,1,-1),& & e_2^{(68)} = (3,1,0,-1),\quad e_3^{(68)} = (9,-1,-1,-1),\, \\ 
e_1^{(83)} = (1,1,1,-1),& & e_2^{(83)} = (3,3,0,-1),\quad e_3^{(83)} = (10,1,-1,-1),\, \\
e_1^{(92)} = (1,1,1,-1),& & e_2^{(92)} = (3,2,0,-1),\quad e_3^{(92)} = (10,0,-1,-1).   
\end{eqnarray*}
Define polytopes $\Delta_1,\, \Delta_1',\, \Delta_2,\, \Delta_2',\, \tilde{\Delta_3},\, \Delta_3$, and $\Delta_3'$ by 
\begin{eqnarray*}
\Delta_1 & := & {\rm Conv}\{ (-1,0,0),\, (0,0,1),\, (-1,2,0),\, (-1,3,-1),\, (2,-1,0)\}, \\
\Delta_1' & := & \Delta^{(42)}={\rm Conv}\left\{\begin{array}{l} (-1,-1,1),\, (-1,-1,-1),\\ (9,-1,-1),\, (0,2,-1),\, (-1,2,-1)\end{array}\right\}, \\
\Delta_2 & := & \Delta^{(68)}={\rm Conv}\{ (-1,0,0),\, (0,0,1),\, (-2,4,-1),\, (-1,3,-1),\, (2,-1,0)\}, \\
\Delta_2' & := & {\rm Conv}\left\{\begin{array}{l} (-1,-1,1),\, (-1,-1,-1),\\ (9,-1,-1),\, (0,2,-1),\, (-1,1,-1)\end{array}\right\}, \\
\tilde{\Delta}_3 & := & {\rm Conv}\{ (-1,0,0),\, (1,-1,1),\, (-1,2,0),\, (-2,4,-1),\, (2,-1,0)\}, \\
\Delta_3 & := & \Delta^{(83)}={\rm Conv}\{ (-1,0,0),\, (0,0,1),\, (-2,4,-1),\, (2,-1,0),\, (1,-1,1)\}, \\
\Delta_3' & := & {\rm Conv}\left\{\begin{array}{l} (-1,-1,1),\, (-1,-1,-1),\\ (9,-1,-1),\, (3,1,-1),\, (-1,2,-1)\end{array}\right\}. 
\end{eqnarray*}
Since the polar dual polytopes $\Delta_1^*$ and $\Delta_2^*$ of $\Delta_1$ and $\Delta_2$, and $\tilde{\Delta}_3^*$ and $\Delta_3^*$ of $\tilde{\Delta}_3$ and $\Delta_3$ are respectively the convex hulls of vertices $(-1,-1,-1)$,\, $(1,0,0)$,\, $(1,3,9)$,\, $(1,3,-1)$, and $(1,0,-1)$, resp. $(-1,-1,-1)$,\, $(1,1,3)$,\, $(1,3,9)$,\, $(1,3,-1)$, and $(1,0,-1)$, resp. $(-1,-1,-1)$,\, $(1,0,-1)$,\, $(1,3,11)$,\, $(1,3,1)$, and $(1,0,-2)$, resp. $(-1,-1,-1)$,\, $(1,0,-1)$,\, $(1,3,11)$,\, $(1,3,1)$, and $(1,1,-1)$, and the linear maps of $\mathbb{R}^3$ determined by matrices $\left(\begin{smallmatrix} 0 & 0 & -1\\ 0 & -1 & -3\\ -1 & -2 & -5 \end{smallmatrix}\right)$, and $\left(\begin{smallmatrix} 0 & 0 & -1\\ 0 & -1 & -4\\ -1 & -2 & -6 \end{smallmatrix}\right)$ give isomorphisms from $\Delta_1'$ to $\Delta_1^*$ and from $\Delta_2'$ to $\Delta_2^*$, respectively from $\Delta_1'$ to $\tilde{\Delta}_3^*$ and from $\Delta_3'$ to $\Delta_3^*$, the relations $\Delta_1^*\simeq\Delta_1'$ and $\Delta_2^*\simeq\Delta_2'$, and $\tilde{\Delta}_3^*\simeq\Delta_1'$ and $\Delta_3^*\simeq\Delta_3'$ hold. 

\noindent
{\bf No. 38.}
The Newton polytope of $F$,which is the convex hull of vertices $(-1,0,0)$,\, $(0,0,1)$,\, $(-1,3,-1)$, and $(2,-1,0)$, is not reflexive. 
It is observed that one might take a polytope with a vertex $(-2,4,-1)$ or $(-1,2,0)$ instead of a face spanned by vertices $(-1,0,0)$,\, $(0,0,1)$, and $(-1,3,-1)$. 
Besides, the Newton polytope of $F'$, which is the convex hull of vertices $(-1,-1,1)$,\, $(-1,-1,-1)$,\, $(4,-1,0)$, and $(0,2,-1)$, is not reflexive. 
It is observed that there are two polytope-dual pairs, that is, $(\Delta_1,\, \Delta^{(42)}=\Delta_1')$, and $(\Delta^{(68)} = \Delta_2,\, \Delta_2')$. 

\noindent
{\bf No. 39.}
The Newton polytope of $F$, which is the convex hull of vertices $(-1,0,0)$,\, $(1,-3,1)$,\, $(-2,4,-1)$, and $(2,-1,0)$, is easily seen to be reflexive. 
Besides, the Newton polytope of $F'$, which is the convex hull of vertices $(-1,-1,1)$,\, $(0,-1,-1)$,\, $(6,0,-1)$, and $(-1,2,-1)$, is not reflexive. 
It is observed that for any reflexive polytope $\Delta$ such that $\Delta_F\subset\Delta\subset\Delta^{(92)}$, the polar dual polytope $\Delta^*$ should be a tetrahedron of which each edge contains $3$ lattice points. 
However, there does not exist a reflexive polytope $\Delta'$ such that $\Delta_{F'}\subset\Delta'\subset\Delta^{(42)}$ that contains a vertex which is adjacent to three vertices between which the edges contain $3$ lattice points. 
Thus, in this case, no pair is polytope-dual. 

\noindent
{\bf No. 40.}
The Newton polytope of $F$, which is the convex hull of vertices $(-1,0,0)$,\, $(1,-1,1)$,\, $(-2,4,-1)$, and $(2,-1,0)$, is not reflexive. 
It is observed that one might take a polytope with a vertex $(0,0,1)$ or $(-1,2,0)$ instead of a face spanned by $(-1,0,0)$,\, $(1,-1,1)$, and $(-2,4,-1)$. 
Besides, the Newton polytope of $F'$, which is the convex hull of vertices $(-1,-1,1)$,\, $(0,-1,-1)$,\, $(9,-1,-1)$, and $(-1,2,-1)$, is not reflexive. 
It is observed that ther are two polytope-dual pairs, that is, $(\tilde{\Delta}_1,\, \Delta^{(42)}=\Delta_1')$, and $(\Delta^{(83)}=\Delta_3,\, \Delta_3')$. 

\subsection{No. 41--No. 43} 
We claim that there exist two polytope-dual pairs for Nos. 41 and 43, and that a unique pair for No. 42. 
Take a basis $\{ e_1^{(n)},\, e_2^{(n)},\, e_3^{(n)}\}$ of lattices $M_{n}$ with  for $n=25,43,48$, and $88$ by 
\begin{eqnarray*}
e_1^{(25)} = (-1,1,0,0), & & e_2^{(25)} = (-3,0,1,0),\quad e_3^{(25)} = (-4,0,0,1),\, \\
e_1^{(43)} = (-1,8,-1,-1),& & e_2^{(43)} = (0,-1,2,-1),\quad e_3^{(43)} = (-1,-1,-1,1),\, \\ 
e_1^{(48)} = (0,8,-1,-1),& & e_2^{(48)} = (-1,-1,2,-1),\quad e_3^{(48)} = (-1,-1,-1,1),\, \\
e_1^{(88)} = (0,4,-1,-1),& & e_2^{(88)} = (-1,-1,2,-1),\quad e_3^{(88)} = (-1,0,-1,1).   
\end{eqnarray*}
Define polytopes $\Delta_1,\, \Delta_1',\, \Delta_2,$ and $\Delta_2'$ by 
\begin{eqnarray*}
\Delta_1 & := & {\rm Conv}\{ (1,0,0),\, (0,1,0),\, (0,0,1),\, (-2,-6,-9)\}, \\
\Delta_1' & := & {\rm Conv}\{ (-1,-1,1),\, (-1,-1,-1),\, (8,-1,-1),\, (-1,2,-1)\}, \\
\Delta_2 & := & {\rm Conv}\{ (1,0,0),\, (0,1,0),\, (0,0,1),\, (0,-2,-3),\, (-1,-3,-4)\}, \\
\Delta_2' & := & \Delta^{(25)} = {\rm Conv}\left\{\begin{array}{l} (-1,-1,1),\, (-1,-1,-1),\\ (8,-1,-1),\, (0,-1,1),\, (-1,2,-1)\end{array}\right\}. 
\end{eqnarray*}
Since the polar dual polytopes $\Delta_i^*$ for $i=1,2$ are the convex hulls of vertices $(-1,-1,1)$,\, $(-1,-1,-1)$,\, $(8,-1,-1)$, and $(-1,2,-1)$, respectively, $(-1,-1,1)$,\, $(-1,-1,-1)$,\, $(8,-1,-1)$,\, $(0,-1,1)$, and $(-1,2,-1)$, the relations $\Delta_1^*=\Delta_1'$ and $\Delta_2^*=\Delta_2'$ clearly hold. 

\noindent
{\bf No. 41.}
The Newton polytope of $F$, which is the convex hull of vertices $(1,0,0)$,\, $(0,1,0)$,\, $(0,0,1)$, and $(-1,-3,-4)$, is not reflexive. 
It is observed that one might take a polytope with a vertex $(-2,-6,-9)$ or $(0,-2,-3)$ instead of a face spanned by vertices $(1,0,0)$,\, $(0,1,0)$, and $(-1,-3,-4)$. 
Besides, the Newton polytope of $F'$, which is the convex hull of vertices $(-1,-1,1)$,\, $(-1,0,-1)$,\, $(8,-1,-1)$, and $(-1,2,-1)$, is not reflexive. 
It is observed that there are two polytope-dual pairs, that is, $(\Delta^{(43)}=\Delta_1,\, \Delta_1')$, and $(\Delta_2,\, \Delta_2')$. 

\noindent
{\bf No. 42.}
The Newton polytope of $F$, which is the convex hull of vertices $(1,0,0)$,\, $(0,1,0)$,\, $(0,0,1)$, and $(-1,-3,-4)$, is not reflexive. 
Besides, the Newton polytope of $F'$, which is the convex hull of vertices $(-1,-1,1)$,\, $(0,-1,-1)$,\, $(4,-1,0)$, and $(-1,2,-1)$, is not reflexive.
It is observed that the pair $(\Delta_2,\, \Delta_2')$ is polytope-dual. 

\noindent
{\bf No. 43.}
The Newton polytope of $F$, which is the convex hull of vertices $(1,0,0)$,\, $(0,1,0)$,\, $(0,0,1)$, and $(-1,-3,-4)$, is not reflexive.
Besides, the Newton polytope of $F'$, which is the convex hull of vertices $(-1,-1,1)$,\, $(0,-1,-1)$,\, $(8,-1,-1)$, and $(-1,2,-1)$, is not reflexive.
It is observed that there are two polytope-dual pairs, that is, $(\Delta^{(48)}=\Delta_1,\, \Delta_1')$, and $(\Delta_2,\, \Delta_2')$. 

\subsection{No. 44} 
We claim that there exists a unique polytope-dual pair. 
Take a basis $\{ e_1^{(n)},\, e_2^{(n)},\, e_3^{(n)}\}$ of lattices $M_{n}$ for $n=7$ and $64$ by 
\begin{eqnarray*}
e_1^{(7)} =  (-1,1,0,0), & & e_2^{(7)} = (-2,0,1,0),\quad e_3^{(7)} = (-4,0,0,1),\, \\
e_1^{(64)} =  (-1,5,-1,-1), & & e_2^{(64)} = (0,-1,2,-1),\quad e_3^{(64)} = (-1,-1,1,0).   
\end{eqnarray*}

Define polytopes $\Delta_1$ and $\Delta_1'$ by 
\begin{eqnarray*}
\Delta_1 & := & {\rm Conv}\{ (1,0,0),\, (0,-1,1),\, (0,0,1),\, (-1,2,-6)\}, \\
\Delta_1' & := & \Delta^{(7)} = {\rm Conv}\{ (-1,-1,1),\, (-1,-1,-1),\, (7,-1,-1),\, (-1,3,-1)\}. 
\end{eqnarray*}
Since the polar dual polytope $\Delta_1^*$ is the convex hull of vertices $(-1,2,1)$,\, $(-1,0,-1)$,\, $(7,0,-1)$, and $(-1,-4,-1)$, and the linear map of $\mathbb{R}^3$ determined by a matrix $\left(\begin{smallmatrix} 1 & 0 & 0\\ 0 & -1 & 0\\ 0 & 1 & 1 \end{smallmatrix}\right)$ gives an isomorphism from $\Delta_1'$ to $\Delta_1^*$, the relation $\Delta_1^*=\Delta_1'$ holds. 

The Newton polytope of $F$, which is the convex hull of vertices $(1,0,0)$,\, $(0,-1,1)$,\, $(0,0,1)$, and $(-1,2,-6)$, coincides with $\Delta_1$. 
Besides, the Newton polytope of $F'$, which is the convex hull of vertices $(-1,-1,1)$,\, $(-1,-1,-1)$,\, $(5,0,-1)$, and $(-1,1,0)$, is not reflexive. 
It is observed that one might take a polytope with vertices $(7,-1,-1)$ and $(-1,b',-1)$ with $b'$ being $0$ or $3$ instead of faces spanned by $(-1,-1,1)$,\, $(-1,-1,-1)$,\, $(5,0,-1)$ and $(-1,1,0)$,\, $(-1,-1,-1)$,\, $(5,0,-1)$. 
It is easily seen that if $b'=3$, there does not exist a polytope in $\Delta^{(64)}$ that contains two pentagonal faces. 
Therefore, the pair $(\Delta_F=\Delta_1,\, \Delta_1')$ is polytope-dual. 

\subsection{No. 45} 
We claim that there exists a unique polytope-dual pair. 
Take a basis $\{ e_1^{(n)},\, e_2^{(n)},\, e_3^{(n)}\}$ of lattices $M_{n}$ for $n=35$ and $66$ by 
\begin{eqnarray*}
e_1^{(35)} = (7,0,-1,-1), & & e_2^{(35)} = (6,-1,0,-1),\quad e_3^{(35)} = (-1,-1,-1,1),\, \\
e_1^{(66)}  = (-1,1,0,0),& & e_2^{(66)} = (-2,0,1,0),\quad e_3^{(66)} = (-3,0,0,1).   
\end{eqnarray*}
Define polytopes $\Delta_1$ and $\Delta_1'$ by 
\begin{eqnarray*}
\Delta_1 & := & {\rm Conv}\{ (0,0,1),\, (0,1,0),\, (2,-2,-1),\, (3,-4,-2),\, (-2,2,-1)\}, \\
\Delta_1' & := & {\rm Conv}\left\{\begin{array}{l} (-1,1,0),\, (0,-1,1),\\ (-1,-1,1),\, (-1,-1,-1),\, (6,-1,-1)\end{array}\right\}. 
\end{eqnarray*}
Since the polar dual polytope $\Delta_1^*$ is the convex hull of vertices $(-1,-1,1)$,\, $(0,-1,-1)$,\, $(-2,-1,-1)$,\, $(-1,0,-1)$, and $(7,6,-1)$, and the linear map of $\mathbb{R}^3$ determined by a matrix $\left(\begin{smallmatrix} -1 & -1 & 0\\ -2 & -2 & 1\\ -4 & -3 & 0 \end{smallmatrix}\right)$ gives an isomorphism from $\Delta_1'$ to $\Delta_1^*$, the relation $\Delta_1^*=\Delta_1'$ holds. 

The Newton polytope of $F$, which is the convex hull of vertices $(0,0,1)$,\, $(0,1,0)$,\, $(3,-4,-2)$, and $(-1,1,0)$, is not reflexive. 
It is observed that one might take a polytope with vertices $(-2,2,-1)$, and $(1,0,0)$ instead of faces spanned by vertices $(0,0,1)$,\, $(0,1,0)$,\, $(3,-4,-2)$, and $(-1,1,0)$,\, $(0,1,0)$,\, $(3,-4,-2)$. 
Besides, the Newton polytope of $F'$, which is the convex hull of vertices $(-1,-1,1)$,\, $(-1,-1,-1)$,\, $(6,-1,-1)$, and $(-1,1,0)$, is not reflexive. 
It is observed that one might take a polytope with a vertex $(0,-1,1)$ instead of a face spanned by vertices $(-1,-1,1)$,\, $(6,-1,-1)$, and $(-1,1,0)$. 
Therefore, the pair $(\Delta_1,\, \Delta_1')$ is polytope-dual. 

\subsection{No. 46--No. 47} 
We claim that there exists a unique polytope-dual pair for Nos. 46 and 47. 
Take a basis $\{ e_1^{(n)},\, e_2^{(n)},\, e_3^{(n)}\}$ of lattices $M_{n}$ for $n=21, 30$ and $86$ by 
\begin{eqnarray*}
e_1^{(21)} =  (-1,1,0,0), & & e_2^{(21)} = (-1,0,1,0),\quad e_3^{(21)} = (-2,0,0,1),\, \\
e_1^{(30)} = (0,4,-1,-1),& & e_2^{(30)} = (-1,-1,4,-1),\quad e_3^{(30)} = (-1,-1,-1,1),\, \\
e_1^{(86)} = (4,0,-1,-1), & & e_2^{(86)} = (3,-1,-1,0),\quad e_3^{(86)} = (0,-1,2,-1).   
\end{eqnarray*}

Define polytopes $\Delta_1,\, \Delta_1',\, \Delta_2$ and $\Delta_2'$ by 
\begin{eqnarray*}
\Delta_1 & := & {\rm Conv}\{ (1,0,0),\, (0,1,0),\, (0,0,1),\, (-1,1,0),\, (2,-3,-1)\}, \\
\Delta_1' & := & {\rm Conv}\left\{ \begin{array}{l}(-1,-1,1),\, (-1,-1,-1),\, (4,-1,-1),\\ (-1,4,-1),\, (-1,2,0),\, (2,-1,0)\end{array}\right\}, \\
\Delta_2 & := & {\rm Conv}\left\{ (1,0,0),\, (0,1,0),\, (0,0,1),\, (-2,-2,-5)\right\}, \\
\Delta_2' & := & {\rm Conv}\{ (-1,-1,1),\, (-1,-1,-1),\, (4,-1,-1),\, (-1,4,-1)\}. 
\end{eqnarray*}
Since the polar dual polytopes $\Delta_1^*$ and $\Delta_2^*$ are the convex hulls of vertices $(0,-1,-1)$,\, $(-1,-1,-1)$,\, $(-1,-1,2)$,\, $(0,-1,4)$,\, $(5,4,-1)$, and $(-1,0,-1)$, respectively, $(-1,-1,-1)$,\, $(-1,-1,1)$,\, $(4,-1,-1)$, and $(-1,4,-1)$, and the linear map of $\mathbb{R}^3$ determined by a matrix $\left(\begin{smallmatrix} -1 & -1 & 1\\ -1 & -1 & 0\\ -3 & -2 & 0 \end{smallmatrix}\right)$ gives an isomorphism from $\Delta_1'$ to $\Delta_1^*$, the relations $\Delta_1^*\simeq\Delta_1'$ and $\Delta_2^*=\Delta_2'$ hold. 

\noindent
{\bf No. 46.}
The Newton polytope of $F$, which is the convex hull of vertices $(0,0,1)$,\, $(0,1,0)$,\, $(2,-3,-1)$, and $(-1,1,0)$, is not reflexive.
It is observed that one might take a polytope with a vertex $(1,0,0)$ instead of a face spanned by vertices $(0,0,1)$,\, $(0,1,0)$, and $(2,-3,-1)$. 
Besides, the Newton polytope of $F'$, which is the convex hull of vertices $(-1,-1,1)$,\, $(-1,0,-1)$,\, $(4,-1,-1)$, and $(-1,2,0)$, is not reflexive. 
It is observed that the pair $(\Delta^{(86)}=\Delta_1,\, \Delta_1')$ is polytope-dual. 

\noindent
{\bf No. 47.}
The Newton polytope of $F$, which is the convex hull of vertices $(1,0,0)$,\, $(0,1,0)$,\, $(0,0,1)$, and $(-1,-1,-2)$, is not reflexive. 
It is observed that one might take a polytope with a vertex $(-2,-2,-5)$ instead of a face spanned by vertices $(1,0,0)$,\, $(0,1,0)$, and $(-1,-1,-2)$. 
Besides, the Newton polytope of $F'$, which is the convex hull of vertices $(-1,-1,1)$,\, $(0,-1,-1)$,\, $(4,-1,-1)$, and $(-1,4,-1)$, is not reflexive. 
It is observed that the pair $(\Delta^{(30)}=\Delta_2,\, \Delta_2')$ is polytope-dual. 

\subsection{No. 48--No. 49} 
We claim that there exists a unique polytope-dual pair. 
Take a basis $\{ e_1^{(n)},\, e_2^{(n)},\, e_3^{(n)}\}$ of lattices $M_{n}$ for $n=5, 56$ and $73$ by 
\begin{eqnarray*}
e_1^{(5)}  =  (-1,1,0,0),& & e_2^{(5)} = (-1,0,1,0),\quad e_3^{(5)} = (-3,0,0,1),\, \\
e_1^{(56)}  =  (-1,0,2,-1), & & e_2^{(56)} = (-1,-1,0,1),\quad e_3^{(56)} = (5,-1,-1,-1),\, \\
e_1^{(73)} = (-1,-1,4,-1), & & e_2^{(73)} = (-1,-1,-1,1),\quad e_3^{(73)} = (5,0,-1,-1).   
\end{eqnarray*}

In \cite{KM}, it is proved that the polytopes $\Delta^{(56)}$ and $\Delta^{(73)}$ are isomorphic to the polytope $\Delta_{(56,73)}$ with vertices $(1,0,0)$,\, $(0,1,0)$,\, $(0,0,1)$, and $(-1,-3,-1)$, under the above choice of basis. 
Let $\Delta_1$ be $\Delta_{(56,73)}$, and $\Delta_1'$ be the polytope $\Delta^{(5)}$ which is the convex hull of vertices $ (-1,-1,1)$,\, $(-1,-1,-1)$,\, $(5,-1,-1)$, and $(-1,5,-1)$. 
Since the polar dual polytope $\Delta_1^*$ of $\Delta_1$ is the convex hull of vertices $(-1,-1,-1)$,\, $(-1,1,-1)$,\, $(5,-1,-1)$, and $(-1,-1,5)$, and the linear map of $\mathbb{R}^3$ determined by a matrix $\left(\begin{smallmatrix} 1 & 0 & 0\\ 0 & 0 & 1\\ 0 & 1 & 0 \end{smallmatrix}\right)$ gives an isomorphism from $\Delta_1'$ to $\Delta_1^*$, the relation $\Delta_1^*\simeq\Delta_1'$ holds. 

\noindent
{\bf No. 48.}
The Newton polytope of $F$, which is the convex hull of vertices $(1,0,0)$,\, $(0,1,0)$,\, $(0,0,1)$, and $(-1,-3,-1)$, coincides with $\Delta_1$. 
Besides, the Newton polytope of $F'$, which is the convex hull of vertices $(-1,-1,1)$,\, $(-1,-1,-1)$,\, $(4,0,-1)$, and $(-1,2,0)$, is not reflexive. 
It is observed that the pair $(\Delta_{F}=\Delta_1,\, \Delta_1')$ is polytope-dual. 

\noindent
{\bf No. 49.}
The Newton polytope of $F$, which is the convex hull of vertices $(1,0,0)$,\, $(0,1,0)$,\, $(0,0,1)$, and $(-1,-3,-1)$, coincides with $\Delta_1$. 
Besides, the Newton polytope of $F'$, which is the convex hull of vertices $(-1,-1,1)$,\, $(-1,-1,-1)$,\, $(4,-1,-1)$, and $(0,4,-1)$, is not reflexive. 
It is observed that the pair $(\Delta_{F}=\Delta_1,\, \Delta_1')$ is polytope-dual. 

\subsection{No. 50} 
We claim that there exists a unique polytope-dual pair. 
Take a basis $\{ e_1^{(n)},\, e_2^{(n)},\, e_3^{(n)}\}$ of lattices $M_{n}$ for $n=1$ and $52$ by 
\begin{eqnarray*}
e_1^{(1)} =  (-1,1,0,0), & & e_2^{(1)} = (-1,0,1,0),\quad e_3^{(1)} = (-1,0,0,1),\, \\
e_1^{(52)} = (-1,2,-1,0), & & e_2^{(52)} = (-1,-1,3,-1),\quad e_3^{(52)} = (-1,-1,-1,2).   
\end{eqnarray*}
Define polytopes $\Delta_1$ and $\Delta_1'$ by 
\begin{eqnarray*}
\Delta_1 & := & {\rm Conv}\{ (1,0,0),\, (0,1,0),\, (0,0,1),\, (-1,-1,-1)\}, \\
\Delta_1' & := & {\rm Conv}\{ (-1,-1,3),\, (-1,-1,-1),\, (3,-1,-1),\, (-1,3,-1)\}. 
\end{eqnarray*}
Since the polar dual polytope $\Delta_1^*$ of $\Delta_1$ is the convex hull of vertices $(-1,-1,3)$,\, $(-1,-1,-1)$,\, $(3,-1,-1)$, and $(-1,3,-1)$, the relation $\Delta_1^*=\Delta_1'$ clearly holds. 

The Newton polytope of $F$, which is the convex hull of vertices $(1,0,0)$,\, $(0,1,0)$,\, $(0,0,1)$, and $(-1,-1,-1)$, coincides with $\Delta_1$. 
Besides, the Newton polytope of $F'$, which is the convex hull of vertices $(0,-1,2)$,\, $(-1,-1,-1)$,\, $(2,-1,-1)$, and $(-1,3,-1)$, is not reflexive. 
It is observed that the pair $(\Delta_{F}=\Delta_1,\, \Delta^{(1)}=\Delta_1')$ is polytope-dual. 

\subsection{No. 51} 
We claim that there exists a unique polytope-dual pair. 
Take a basis $\{ e_1^{(n)},\, e_2^{(n)},\, e_3^{(n)}\}$ of the lattice $M_{32}$ by 
\begin{eqnarray*}
e_1^{(32)} = (-1,0,3,-1), & &  e_2^{(32)} = (0,-1,3,-1),\quad e_3^{(32)} = (-1,-1,-1,1).   
\end{eqnarray*}
Define polytopes $\Delta_1$ and $\Delta_1'$ by 
\begin{eqnarray*}
\Delta_1 & := & {\rm Conv}\{ (0,0,1),\, (1,0,0),\, (-2,2,-1),\, (-4,3,-2),\, (2,-2,-1)\}, \\
\Delta_1' & := & {\rm Conv}\{ (0,0,1),\, (1,0,0),\, (0,1,0),\, (-4,3,-2),\, (3,-4,-2)\}. 
\end{eqnarray*}
Since the polar dual polytope $\Delta_1^*$ of $\Delta_1$ is the convex hull of vertices $(-1,-1,1)$,\, $(-1,0,-1)$,\, $(-1,-2,-1)$,\, $(0,-1,-1)$, and $(6,7,-1)$, and the linear map of $\mathbb{R}^3$ determined by a matrix $\left(\begin{smallmatrix} -1 & -2 & -1\\ 0 & -1 & -1\\ -1 & -1 & 1 \end{smallmatrix}\right)$ gives an isomorphism from $\Delta_1'$ to $\Delta_1^*$, the relation $\Delta_1^*\simeq\Delta_1'$ holds. 

The Newton polytope of $F=F'$, which is the convex hull of vertices $(1,0,0)$,\, $(-4,3,-2)$,\, $(0,0,1)$, and $(2,-2,-1)$, is not reflexive. 
It is observed that one might take a polytope with a vertex $(0,1,0)$ or $(-2,2,-1)$ instead of a face spanned by vertices $(1,0,0)$,\, $(0,0,1)$, and $(-4,3,-2)$. 
Therefore there the pair $(\Delta_1,\, \Delta^{(32)}=\Delta_1')$ is polytope-dual. 

Therefore, the claims are verified. 
\QED

\section{Closing Remarks}
As is mentioned in~\cite{Ebeling}, almost all mirror symmetric pairs of weight systems  in the sense of~\cite{Belcastro} are also strongly coupled. 
We would like to study not only the full families of $K3$ surfaces, but families of $K3$ surfaces associated to reflexive polytopes $\Delta$ and $\Delta'$ obtained in Theorem~\ref{MainThm}. 
More precisely, it is interesting to study a relation with a lattice duality~\cite{DolgachevMirror} and the coupling.

Makiko Mase\\
Universit\"at Mannheim, Lehrstuhl f\"ur Mathematik VI\\
B6, 26, 68131 Mannheim, Germany \\
email: mmase@mail.uni-mannheim.de
\end{document}